\documentclass[11pt]{article}
\usepackage[margin=1in]{geometry}
\usepackage[utf8]{inputenc}
\usepackage{amsmath,amssymb,amsthm}
\usepackage{color}
\usepackage{graphicx}
\usepackage{caption}
\usepackage{subcaption}

\usepackage{lipsum}
\newcommand\blfootnote[1]{%
  \begingroup
  \renewcommand\thefootnote{}\footnote{#1}%
  \addtocounter{footnote}{-1}%
  \endgroup
}

\newtheorem{theorem}{Theorem}
\newtheorem{corollary}[theorem]{Corollary}
\newtheorem{lemma}[theorem]{Lemma}
\newtheorem{proposition}[theorem]{Proposition}
\newtheorem{remark}{Remark}

\title{Data-Driven Moment-Based Distributionally Robust Chance-Constrained Optimization}
\author{Joshua Comden \hspace{0.5in}
Ahmed S. Zamzam \hspace{0.5in}
Andrey Bernstein}
\date{}

\begin{document}

\maketitle

\begin{abstract}
    Many stochastic optimization problems include chance constraints that enforce constraint satisfaction with a specific probability; 
however, solving an optimization problem with chance constraints assumes that the solver has access to the exact underlying probability distribution, which is often unreasonable.
In data-driven applications, it is common instead to use historical data samples as a surrogate to the distribution;
however, this comes at a significant computational cost from the added time spent either processing the data or, worse, adding additional variables and constraints to the optimization problem.
On the other hand, the sample mean and covariance matrix are lightweight to calculate, and it is possible to reframe the chance constraint as a distributionally robust chance constraint.
The challenge here is that the sample mean and covariance matrix themselves are random variables, so their uncertainty should be factored into the chance constraint.
This work bridges this gap by modifying the standard method of distributionally robust chance constraints to guarantee its satisfaction.
The proposed data-driven method is tested on a particularly problematic example. The results show that the computationally fast proposed method is not significantly more conservative than other methods. 
\end{abstract}

\blfootnote{J. Comden, A. Zamzam, and A. Bernstein are with the National Renewable Energy Laboratory (NREL), Golden CO, USA (Email: \{joshua.comden, ahmed.zamzam, andrey.bernstein\}@nrel.gov).}

\blfootnote{The views expressed in the article do not necessarily represent the views of the DOE or the U.S. Government. The U.S. Government retains and the publisher, by accepting the article for publication, acknowledges that the U.S. Government retains a nonexclusive, paid-up, irrevocable, worldwide license to publish or reproduce the published form of this work, or allow others to do so, for U.S. Government purposes.}

\section{Introduction}
\label{sec:intro}
Many stochastic optimization problems include probabilistic constraints, also called chance constraints. Such constraints are used to restrict the feasible region of the decision variables such that the probability of a certain stochastic constraint being satisfied is high enough.
The main motivation for using a chance constraint is to allow a regulated extent of constraint violation; without an allowable violation, the objective function might be significantly degraded or the decision space might not be feasible for all stochastic realizations. 
Chance constraints also play an important role in sequential decision problems where an optimization problem must be solved at each time step---e.g., chemical process control \cite{schwarm1999chance}, water network management \cite{grosso2014chance}, spacecraft rendezvous \cite{gavilan2012chance}, pharmacy inventory management \cite{jurado2016stock}, and energy management \cite{darivianakis2017power}, among many others.
Thus, it is vital that chance-constrained optimization problems are solved within a reasonable amount of time in sequential optimization settings.

Much of the earlier research on chance constraints assumes that the probability distributions are known beforehand, and hence, methods have been developed for particular types of distributions and applications.
For instance, if the random variables in an application are assumed to be Gaussian with known moments, then the chance constraint lends itself to a simple deterministic second-order cone reformulation  (e.g., \cite{summers2015stochastic,venzke2017convex}).
In some cases, the exact shape of the probability distribution might not be known, but there could be some properties that are assumed to be known (e.g., its modality and support).
Thus, it is desirable to guarantee that the chance constraint is satisfied for any probability distribution within a family containing all distributions with the assumed properties.
This is known as a distributionally robust chance constraint, or an ambiguous chance constraint, which is defined by the family of distributions for which the chance constraint needs to be satisfied.

\subsection{Data-Driven Chance-Constrained Optimization}
For many applications, only an empirical distribution from historical data is obtainable. This is more evident in sequential decision problems that use system measurements to make decisions.
Recently, there has been significant interest in tackling chance-constrained problems when only samples from the random variables are available \cite{calafiore2006distributionally,stellato2014data,hong2017learning,yanikouglu2013safe,delage2010distributionally,bertsimas2018data,roveto2020co,poolla2020wasserstein,coulson2020distributionally,erdougan2006ambiguous}.
For general mutually independent random variables, \cite{nemirovski2007convex} use the data to build conservative approximations of the chance constraint(s).
If the random variables are assumed to be Gaussian and mutually independent, \cite{darivianakis2017power,mieth2018data} use confidence intervals around the sample means and variances to build a deterministic approximation of the chance constraint.
The work by \cite{bertsimas2018data} adapts different statistical hypothesis tests to convert data and a chance constraint into deterministic constraints.
For instance, they use the goodness-of-fit hypothesis test from \cite{bertsimas2018robust} to add deterministic constraints, which unfortunately grow linearly with the number of samples.
The data-driven method developed by \cite{hong2017learning} uses the data to build a convex uncertainty set to replace the chance constraint in the optimization problem.
The uncertainty set is built by combining polytopes and ellipsoids, which can be as complex or simple as the optimization practitioner decides.
The authors in \cite{jiang2016data} show that a distributionally robust chance constraint defined by any general $\phi$-divergence from the empirical distribution is equivalent to a standard chance constraint defined by the empirical distribution with a perturbed bound on the probability of violation.
That allows the practitioner to convert a $\phi$-divergence tolerance from the distributionally robust chance constraint into a perturbation on the bound for a standard chance constraint that can then be solved with regular chance-constrained methods;
however, the practitioner now has the flexibility (or burden) to decide the value of the tolerance, even with the method's given designed guidelines.

The most common and simple distributionally robust chance constraints are those defined by a family of all distributions with a given mean and covariance.
The main reason for being popular is that if the true mean and covariance are known, then the distributionally robust chance constraint has a simple convex deterministic equivalent that can be easily coded into an optimization problem solver (See Lemma \ref{thm:calafiore3p1all} in Section \ref{sec:dist_rob_cc} by \cite{ghaoui2003worst} and \cite{calafiore2006distributionally}).
In fact, this specific form lends itself to an elegant way to price risk and variability \cite{mieth2020risk}.
Because the true mean and covariance are generally not known, a practical work-around is to simply use the sample mean and covariance from historical data in place of the true ones (e.g., \cite{baker2016distribution,xie2017distributionally});
however, this means that the original distributionally robust chance constraint is not necessarily guaranteed to hold because the sample mean and covariance are random variables themselves.
If the true mean and covariance are instead assumed to be contained within a convex set around the sampled versions, then the chance constraint can be converted into a semidefinite program~\cite{ghaoui2003worst}.
Our proposed method is similar to the data-driven method developed by \cite{calafiore2006distributionally} in that it uses the data to calculate the sample mean, sample covariance, and confidence regions around them to form a deterministic convex surrogate constraint.
With the same general idea, the method from \cite{bertsimas2018data} modifies the confidence region designed by \cite{delage2010distributionally} to form a deterministic convex constraint.

One particular downside to all the previously described data-driven methods is that they require the optimization practitioner to choose a probability at which the chance constraint itself can be violated.
In other words, the true probability of violation for the original stochastic constraint might be larger than the bound set by the chance constraint.
Also, when the chance constraint fails, there is not necessarily a bound on the violation of the chance constraint.
Thus, the practitioner would be required to tolerate this increased probability of violation or need to make adjustments to the chance constraint to counteract it.
For these methods, however, it is unclear what probability the practitioner should allow for chance-constraint violation or how to adjust the originally intended chance constraint to offset this violation probability.
On the other hand, the method by \cite{stellato2014data}, called the Multivariate Sampled Chebyshev Approach, is moment-based and does not require such a probability to be chosen.
The downside to this method is that the minimum number of samples needed increases linearly with the dimension of the random vector, and it does not converge to the deterministic equivalent  when the moments are known a priori.
From the opposite direction of allowing a probability of chance-constraint violation, the method developed by \cite{bienstock2014chance}, which assumes Gaussian independent random variables, first solves the optimization problem as if the sample means and variances are the true values.
Afterward, they use out-of-sample analysis to calculate the true probability of violation of the stochastic constraint.

\subsection{Computational Complexity Advantage of Moment-Based Methods}

With the recent use of larger data sets to represent empirical distributions, the added computational time needed to find a solution becomes more important, especially for sequential decision problems that require a solution in a timely matter.
Data-driven chance-constrained optimization methods either process the data before solving the optimization problem or integrate the data directly into the optimization problem by introducing additional decision variables and constraints.
Although the methods that do the latter have the potential to yield less conservative solutions, the added computational time might render them impractical.

We demonstrate this fact with a simple sequential decision problem and compare the computational complexity of the different methods.
Suppose that an application requires an $n-$dimensional decision to be made at each time step from a linear program that includes a chance constraint.
The chance constraint is dependant on an $m-$dimensional vector of sampled random variables, where $N$ is the number of samples so far.
Table \ref{tab:cc_complexity} gives the computational complexity for processing a single new sample and solving the optimization problem for different methods.
Only the computational time of the moment-based, moment uncertainty set and the union of discretized cells methods do not increase with increasing sample size;
however, the latter requires a second-order cone program to be solved $O(\sqrt{m})$ times.
Moment-based methods require only $\theta(m)$ time to update the sample mean and $\theta(m^2)$ to update the sample covariance with the new sample while it replaces the chance constraint with a second-order cone constraint.
Thus, moment-based methods that require only the sample mean (e.g., the method described in Section \ref{sec:indep_data-driven_guar}) require only $\theta(m)$ time to process a new sample.
The sample processing time of the uncertainty set method increases with sample size because of an ordering step: either the size of the uncertainty set can be updated by inserting the sample into the previous ordering, or the shape of the uncertainty can be updated, which will require a reordering of all the samples.
All remaining methods add $\Omega(N)$ auxiliary variables and constraints to the optimization problem.
The calculated running times in the table for solving the optimization problem are based on the assumption that a standard linear program can be solved in $O(n^3)$ time \cite{gonzaga1989algorithm,vaidya1990algorithm} and that a standard second-order cone program or semidefinite program can be solved in $O(n^{3.5})$ time. (The algorithms by \cite{monteiro2000polynomial,monteiro1998polynomial} requires $O(\sqrt{n})$ iterations with each iteration requiring an inverse of an $n\times n$ matrix.)

\begin{table}
    \caption{Computational complexity for data-driven chance-constrained methods.\label{tab:cc_complexity}} 
    {\begin{tabular}{|p{0.4\textwidth}||p{0.25\textwidth}|p{0.25\textwidth}|}
        \hline
        \hfil \textbf{Chance-Constrained Method} & \hfil Processing New Sample & \hfil Solving Opt. Problem \\
        \hline \hline
        \hfil Moment-based  & \hfil $\theta(m^2)$ & \hfil $O(n^{3.5})$ \\
        \hline
        \hfil Uncertainty set & \hfil $\Omega(m^2+\log N)$  & \hfil $O(n^{3.5})$ \\
        \hline
        \hfil Union of discretized cells & \hfil $\theta(m)$ & \hfil $O(\sqrt{m}(n+m)^{3.5})$ \\
        \hline
        \hfil Moment uncertainty set & \hfil $O(m^3)$ & \hfil $O((n+m^2)^{3.5})$ \\
        \hline\hline
        \hfil Data-driven CVaR & \hfil $O(m)$ & \hfil $O((n+N)^{3.5})$ \\
        \hline
        \hfil Linear-convex ordering & \hfil --- & \hfil $O((n+m+N)^3)$ \\
        \hline
        \hfil Prohorov metric & \hfil --- & \hfil $O((nN)^3)$ \\
        \hline
    \end{tabular}}
    {\emph{Moment-based} methods include \cite{calafiore2006distributionally}, MSCA by \cite{stellato2014data}, and the method proposed in this paper.  The \emph{uncertainty set} method was developed by \cite{hong2017learning}. The \emph{union of discretized cells} was developed by \cite{yanikouglu2013safe}.  The \emph{moment uncertainty set} method was developed by \cite{delage2010distributionally} and adapted by \cite{bertsimas2018data}.  \emph{Data-driven CVaR} methods include those developed by \cite{roveto2020co}, \cite{poolla2020wasserstein}, and \cite{coulson2020distributionally}.  The \emph{linear-convex ordering} method comes from \cite{bertsimas2018robust} being adapted by \cite{bertsimas2018data}.  The \emph{Prohorov metric} method comes from \cite{erdougan2006ambiguous}. }
\end{table}

\subsection{Contributions}

The developed approaches in this paper aid in the practical application and implementation of chance constraints through the use of available data samples, especially for sequential decision problems. The main contributions can be summarized as:
\begin{enumerate}
    \item We design a moment-based distributionally robust method that can automatically adjust the chance constraint to completely account for the uncertainty from having only a limited number of samples to estimate the mean and covariance.
    \item We prove that the chance constraint is guaranteed to be satisfied by the method and that it asymptotically approaches the distributionally robust chance constraint as if the true mean and covariance were known a priori.
    \item For the special case when the random variables are assumed to be independent, we provide a method that does not require estimating the covariance, and one that does require this estimation. We prove the satisfaction of the chance constraint using the same proof structure for both methods, and we contrast their differences.
    \item We demonstrate the importance of having these guarantees with a numerical evaluation of a sports betting scenario where bets need to be allocated among correlated wagers.  It shows that the proposed method is not significantly more conservative than other methods while giving a theoretical guarantee of chance-constraint satisfaction and reduced computation time.
\end{enumerate}

\section{Distributionally Robust Chance Constraints}
\label{sec:dist_rob_cc}
Let $\mathbf{x}\in\mathbb{R}^n$ be a vector of decision variables, and let $\mathbf{a}\in\mathbb{R}^{n}$ be a vector of random variables drawn from the probability distribution $P^*$.
For the decision variables $\mathbf{x}$, a linear chance constraint can be stated as:
\begin{align}
    \text{Pr}\left(\mathbf{a}^\intercal\mathbf{x}\leq 0|\mathbf{a}\sim P^*\right) \geq 1-\alpha \label{eq:prob_cc_original}
\end{align}
where $\alpha\in(0,1)$ is a scalar that sets an upper limit on the probability of $\mathbf{a}^\intercal\mathbf{x}\leq 0$ not being satisfied, i.e., the probability of constraint violation.


\begin{remark}
    Although the linear inequality constraint $\mathbf{a}^\intercal\mathbf{x}\leq 0$ might seem restrictive, with the use of a vector of auxiliary variables $\mathbf{y}$, it can conservatively approximate an inequality of a sum of randomly weighted functions, i.e., $\sum_{i=1}^n a_ig_i(\mathbf{x}) \leq 0$:
    \begin{align}
        \text{Pr}\left(\mathbf{a}^\intercal\mathbf{y}\leq 0|\mathbf{a}\sim P^*\right) & \geq 1-\alpha \label{eq:gen_con} \\
        g_i(\mathbf{x}) & \leq y_i, \quad \forall i\in\{1,\dots,l\} \nonumber \\
        g_i(\mathbf{x}) & \geq y_i, \quad \forall i\in\{l+1,\dots,l+p\} \nonumber \\
        g_i(\mathbf{x}) & = y_i, \quad \forall i\in\{l+p+1,\dots,n\} \nonumber
    \end{align}
    if the random variables $\{a_1,\dots,a_l\}$ are nonnegative and the random variables $\{a_{l+1},\dots,a_{l+p}\}$ are nonpositive.
    Any function $g_i$ could be replaced with a constant to capture additive randomness independent of the decision $\mathbf{x}$.
    If the inequality constraints for the auxiliary variables are tight---i.e., $g_i(\mathbf{x})=y_i:\forall i\in\{1,\dots,l+p\}$---then it is equivalent to the original chance constraint.
    Additionally, if the functions $\{g_{1}(\mathbf{x}),\dots,g_{l}(\mathbf{x})\}$ are convex, $\{g_{l+1}(\mathbf{x}),\dots,g_{l+p}(\mathbf{x})\}$ are concave, $\{g_{l+p+1}(\mathbf{x}),\dots,g_{n}(\mathbf{x})\}$ are linear, and the linear chance constraint is replaced by its convex equivalent if one exists or a convex approximation, then \eqref{eq:gen_con} is a set of convex constraints.
\end{remark}

Typically, knowing the true probability distribution $P^*$ is not practically valid.
On the other hand, it is more likely for certain properties of $P^*$ to be known or estimated, such as its mean and covariance, and then the chance constraint can be formulated to ensure satisfaction for any distribution with the same properties, leading to 
a distributionally robust chance constraint.
Although many different properties can be used to define a family of distributions, a very common and intuitive set of properties to use is the first two moments.
This is the direction we pursue here, which is often referred to as the Chebyshev ambiguity set~\cite{hanasusanto2015distributionally}.

Let $\mathcal{P}(\boldsymbol{\mu},\boldsymbol{\Sigma})$ be the family of all probability distributions with mean $\boldsymbol{\mu}$, and covariance $\boldsymbol{\Sigma}$.
Later, we will slightly abuse this notation by having a third or fourth argument in $\mathcal{P}$, which will be either the support of the distribution or bounds on the support; this will be specified when it occurs. 
We use $\boldsymbol{\mu}^*$ and $\boldsymbol{\Sigma}^*$ to denote the true mean and covariance of the distribution $P^*$, respectively. Then, the distributionally robust chance constraint defined by the family $\mathcal{P}(\boldsymbol{\mu}^*,\boldsymbol{\Sigma}^*)$ can be formulated as:
\begin{align}
    \inf_{P\in\mathcal{P}(\boldsymbol{\mu}^*,\boldsymbol{\Sigma}^*)}\text{Pr}\left(\mathbf{a}^\intercal\mathbf{x}\leq 0|\mathbf{a}\sim P\right) \geq 1- \alpha  \label{eq:prob_cc_distrob_true}
\end{align}
which conservatively guarantees that the original chance constraint \eqref{eq:prob_cc_original} is satisfied.
If the true mean $\boldsymbol{\mu}^*$ and covariance $\boldsymbol{\Sigma}^*$ are known, the following Lemma from \cite{calafiore2006distributionally} gives a convex deterministic constraint equivalent to \eqref{eq:prob_cc_distrob_true}.

\begin{lemma}
    \label{thm:calafiore3p1all}
    \cite{ghaoui2003worst}(Theorem 1), \cite{calafiore2006distributionally}(Theorem 3.1) The distributionally robust chance constraint \eqref{eq:prob_cc_distrob_true} for $\alpha\in(0,1)$ is equivalent to:
    \begin{align}
        \boldsymbol{\mu}^{*\intercal}\mathbf{x} + \sqrt{\frac{1-\alpha}{\alpha}}\sqrt{\mathbf{x}^\intercal\boldsymbol{\Sigma}^*\mathbf{x}} \leq 0. \label{eq:calafiore3p1all}
    \end{align}
\end{lemma}

In many practical applications, the true mean $\boldsymbol{\mu}^*$ and covariance $\boldsymbol{\Sigma}^*$ are not known but instead can be estimated from data samples.
The next section focuses on guaranteeing the distributionally robust chance constraint \eqref{eq:prob_cc_distrob_true} similar to Lemma \ref{thm:calafiore3p1all} by leveraging only samples from the distribution $P^*$.

\section{Data-Driven Guarantees}
\label{sec:data-driven_gaur}

Let $\mathbf{a}_1,\dots,\mathbf{a}_N$ be $N$ independent samples drawn from the probability distribution $P^*$, which has a support $\mathcal{S}^*$.
From the data samples, we can find the sample mean as $\hat{\boldsymbol{\mu}}_N:=\frac{1}{N}\sum_{i=1}^N\mathbf{a}_i$ and the sample covariance as $\hat{\boldsymbol{\Sigma}}_N:=\frac{1}{N}\sum_{i=1}^N(\mathbf{a}_i-\hat{\boldsymbol{\mu}}_N)(\mathbf{a}_i-\hat{\boldsymbol{\mu}}_N)^\intercal$.

In addition to the data, we assume that we have an estimate $\hat{\mathcal{S}}$ of the true support $\mathcal{S}^*$ such that it contains the true support, i.e., $\mathcal{S}^*\subseteq\hat{\mathcal{S}}$.
This assumption will be relaxed in Section \ref{sec:relax_support}.
We define the radius of the estimated support $\hat{\mathcal{S}}$ as a function of $\mathbf{x}\in\mathbb{R}^n$ that acts as a scaling vector:
\begin{align}
    r(\mathbf{x}):= \frac{1}{2}\sup_{\{\mathbf{a}_1,\mathbf{a}_2\}\in\hat{\mathcal{S}}}|\mathbf{a}_1^\intercal\mathbf{x}-\mathbf{a}_2^\intercal\mathbf{x}| \label{eq:rx_general}
\end{align}
which is a convex function even if $\hat{\mathcal{S}}$ is not convex (See \cite{boyd2004convex}, Example 3.7).
We next show the closed form of \eqref{eq:rx_general}
for a couple common examples of $\hat{\mathcal{S}}$.

If $\hat{\mathcal{S}}$ is a set of independent intervals, i.e., a hyperbox, then each element of $\mathbf{a}_1$ and $\mathbf{a}_2$ can be optimized independently.
Let $\underline{\mathbf{s}}$ and $\overline{\mathbf{s}}$ be the vectors collecting the upper and lower bounds in all dimensions in $\mathcal{\hat{S}}$.
Then, the radius is the following weighted L1-norm of $\mathbf{x}$:
\begin{align}
    r_\text{int}(\mathbf{x}) = \frac{1}{2}\sum_{i=1}^n|x_i|(\overline{s}_i-\underline{s}_i) = \frac{1}{2}\|\mathbf{S}\mathbf{x}\|_1 \label{eq:rx_indint}
\end{align}
where $\mathbf{S}:=\text{diag}(\overline{\mathbf{s}}-\underline{\mathbf{s}})$.
In addition, if $\hat{\mathcal{S}}$ is a polytope and all of its vertices are enumerated in the set $(\mathbf{v}_1,\dots,\mathbf{v}_m)$, then the radius is the difference between the maximum and minimum inner products of the vertices with $\mathbf{x}$:
\begin{align}
    r_\text{poly}(\mathbf{x}) = \frac{1}{2}\max \left\{\mathbf{x}^\intercal\mathbf{v}_1,\dots,\mathbf{x}^\intercal\mathbf{v}_m\right\} -\frac{1}{2}\min \left\{\mathbf{x}^\intercal\mathbf{v}_1,\dots,\mathbf{x}^\intercal\mathbf{v}_m\right\}. \nonumber
\end{align}
If $\hat{\mathcal{S}}$ is the ellipsoid $(\mathbf{a}-\mathbf{c})^\intercal\mathbf{V}(\mathbf{a}-\mathbf{c})\leq 1$, where $\mathbf{c}$ is its center and $\mathbf{V}$ is a symmetric positive definite matrix,
then the radius is given by the following norm:
\begin{align}
    r_\text{ellip}(\mathbf{x}) = \| \mathbf{x} \|_{\mathbf{V}} =  \sqrt{\mathbf{x}^\intercal\mathbf{V}^{-1}\mathbf{x}}. \nonumber
\end{align}
See the Appendix \ref{prf:r_ellip} for the proof.

\subsection{General Guarantees}

With the sample mean $\hat{\boldsymbol{\mu}}_N$, sample covariance $\hat{\boldsymbol{\Sigma}}_N$, and the radius $r(\mathbf{x})$ of the estimated support $\hat{\mathcal{S}}$ as a function of the decision variable $\mathbf{x}$, we replace the distributionally robust chance constraint \eqref{eq:prob_cc_distrob_true} in an optimization problem with the following deterministic convex constraint:
\begin{subequations}
    \label{eq:prob_cc_distrob_data}
    \begin{align}
        \hat{\boldsymbol{\mu}}_N^\intercal\mathbf{x} + \phi_N r(\mathbf{x}) + \kappa_N\sqrt{\frac{1-\alpha}{\alpha}}\|\mathbf{y}\|_2 & \leq 0 \\
        \sqrt{\mathbf{x}^\intercal\hat{\boldsymbol{\Sigma}}_N\mathbf{x}} & \leq y_1 \label{eq:prob_cc_distrob_data_y1} \\
        \sqrt{2\phi_N}r(\mathbf{x}) & \leq y_2 \label{eq:prob_cc_distrob_data_y2}
    \end{align}
\end{subequations}
where $\mathbf{y}\in\mathbb{R}^2$ is a vector of auxiliary variables, and $(\phi_N,\kappa_N)$ are positive scalars that depend on the number of samples $N$ and $\alpha$, to be determined later.
This deterministic constraint is a conservative approximation that guarantees \eqref{eq:prob_cc_distrob_true} will be satisfied and asymptotically approaches it as the number of samples $N$ approaches infinity.
The following theorem proves these statements for specific settings of the scalars $(\phi_N,\kappa_N)$ under a mild condition on the number of samples $N$ and $\alpha$.

\begin{theorem}
    \label{thm:data_to_chance_const}
    Assume that the estimated support $\hat{\mathcal{S}}$ contains the true support $\mathcal{S}^*$, i.e., $\mathcal{S}^*\subseteq\hat{\mathcal{S}}$.
    Let there exist a scalar $p > 2$ such that:
        \begin{align}
            N & >\left(2+\sqrt{2\ln(4/\alpha)}\right)^p \label{eq:proof2_DRCC_conditions_a} 
        \end{align}
    where $\alpha\in(0,1)$.
    Define the following scalars that depend on $N$, $\alpha$, and $p$:
    \begin{subequations}
        \label{eq:proof2_DRCC_constants}
        \begin{align}
            \kappa_N & := \left(1-\frac{4}{\alpha}\exp\left(-(N^\frac{1}{p}-2)^2/2\right)\right)^{-\frac{1}{2}} \\
            \phi_N & := N^{\left(\frac{1}{p}-\frac{1}{2}\right)}.
        \end{align}
    \end{subequations}
    Then, the convex deterministic constraints \eqref{eq:prob_cc_distrob_data}
    represent a conservative approximation of the distributionally robust chance constraint \eqref{eq:prob_cc_distrob_true}. In addition,  the constraints \eqref{eq:prob_cc_distrob_data} asymptotically approach \eqref{eq:prob_cc_distrob_true} as $N\rightarrow\infty$.
\end{theorem}
A key aspect of the proof comes from the following lemmas.
Lemma \ref{thm:confregion} uses a confidence region of the first two moments to guarantee that a univariate version of  the  distributionally robust chance constraint \eqref{eq:prob_cc_distrob_true} will be satisfied.
Then the confidence region described by Lemma \ref{thm:calafiore4p1} is applied.
Finally, the probability that defines the confidence region is set depending on the number of samples $N$ and $\alpha$.
See the Appendix \ref{prf:data_to_chance_const} for the full proof.

\begin{lemma}
    \label{thm:confregion}
    Let $c$ be any scalar and $v$ be a random variable distributed according to a probability distribution in the family $\mathcal{P}(\mu^*,\sigma^{2*},\underline{v}^*,\overline{v}^*)$, where the mean $\mu^*$ and variance $\sigma^{2*}$ are fixed, and the support is $[\underline{v}^*,\overline{v}^*]$.
    Let $\mathcal{Q}_\delta$ be a confidence region of the first two moments such that $\mathcal{Q}_\delta$ contains $(\mu^*,{\sigma^*}^2)$ with a probability of at least $1-\delta$ for $\delta\in(0,1)$.
    If the following is satisfied:
    \begin{align}
        \inf_{P\in\mathcal{P}(\mu,\sigma^2,\underline{v}^*,\overline{v}^*)}\text{Pr}\left(v+c \leq 0 | v\sim P \right) \geq 1-\epsilon, \quad \forall (\mu,\sigma^2)\in \mathcal{Q}_\delta \label{eq:proof2_const_DRO_samp}
    \end{align}
    for $\epsilon\in(0,1)$, then: 
    \begin{align}
        \inf_{P\in\mathcal{P}(\mu^*,\sigma^{2*},\underline{v}^*,\overline{v}^*)}\text{Pr}\left(v+c \leq 0 | v\sim P \right) \geq (1-\epsilon)(1-\delta). \label{eq:proof2_const_DRO_perf}
    \end{align}
\end{lemma}
See the Appendix \ref{prf:confregion} for the proof.

\begin{lemma}
    \label{thm:calafiore4p1}\cite{calafiore2006distributionally}(In proof of Theorem 4.1) Let $\mathbf{v}_1,\dots\mathbf{v}_N$ be $N\geq \left(2+\sqrt{2\ln(4/\delta)}\right)^2$ independent samples drawn from the probability distribution $P^*$ with support $\mathcal{S}^*$ for $\delta\in(0,1)$.
    Define $r^*$ as the radius of $\mathcal{S}^*$, i.e. $r^*:=\frac{1}{2}\sup_{\{\mathbf{v}_1,\mathbf{v}_2\}\in\mathcal{S}^*}\|\mathbf{v}_1-\mathbf{v}_2\|_2$.
    Then,
    \begin{subequations}
        \begin{align}
            \|\boldsymbol{\mu}^*-\hat{\boldsymbol{\mu}}_N\|_2 & \leq \frac{r^*}{\sqrt{N}}\left(2+\sqrt{2\ln(2/\delta)}\right) \\
            \|\boldsymbol{\Sigma}^*-\hat{\boldsymbol{\Sigma}}_N\|_F & \leq \frac{2(r^*)^2}{\sqrt{N}}\left(2+\sqrt{2\ln(4/\delta)}\right)
        \end{align}
    \end{subequations}
    with a probability of at least $1-\delta$.
\end{lemma}

For a practitioner, one downside to Theorem \ref{thm:data_to_chance_const} is that the value of $p$ needs careful tuning.
A larger value of $p$ in Equation \eqref{eq:proof2_DRCC_constants} would result in $\phi_N$ approaching zero faster as the number of samples $N$ increases;
however, this has the opposite effect on $\kappa_N$, which approaches one, and hence it increases the number of samples needed to satisfy the condition \eqref{eq:proof2_DRCC_conditions_a}.
Instead, the following corollary sets $p$ as a function of $N$ and $\alpha$, so $(\kappa_N,\phi_N)$ are completely defined by $N$ and $\alpha$, which relieves the practitioner from needing to decide $p$.
\begin{corollary}
    \label{thm:data_to_chance_const_p1}
    Theorem \ref{thm:data_to_chance_const} remains true if the condition \eqref{eq:proof2_DRCC_conditions_a} is replaced by:
    \begin{align}
        \sqrt{\frac{16N}{\exp{\left(\left(\sqrt{N}-2\right)^2\right)}}} < \alpha \label{eq:proof2_N_alpha_p1}
    \end{align}
    and the scalar definitions in \eqref{eq:proof2_DRCC_constants} are replaced with:
    \begin{subequations}
        \begin{align}
            \kappa_N & := \sqrt{\frac{\sqrt{N}}{\sqrt{N}-1}} \\
            \phi_N & := \frac{2+\sqrt{2\ln (4\sqrt{N}/\alpha)}}{\sqrt{N}}.
        \end{align}
    \end{subequations}
\end{corollary}
\proof
    Let $p:=\log_{2+\sqrt{2\ln (4\sqrt{N}/\alpha)}}(N)$ in Theorem \ref{thm:data_to_chance_const}, then the corollary comes up.
    Specifically, with this setting of $p$, we have that $N^{\frac{1}{p}} = 2+\sqrt{2\ln (4\sqrt{N}/\alpha)}$, and hence we get $\kappa_N=\sqrt{\frac{\sqrt{N}}{\sqrt{N}-1}}$ and $\phi_N := \frac{2+\sqrt{2\ln (4\sqrt{N}/\alpha)}}{\sqrt{N}}$.
    This value of $N^{\frac{1}{p}}$ means that the lower bound condition on $N$ in Theorem \ref{thm:data_to_chance_const} is satisfied for $N\geq 1$.
    Also, the lower bound on $\alpha$ in \eqref{eq:proof2_N_alpha_p1} ensures that $p > 2$, which can be shown by plugging the LHS of \eqref{eq:proof2_N_alpha_p1} in this specific setting of $p$. 
\endproof

Although condition \eqref{eq:proof2_DRCC_conditions_a} in Theorem \ref{thm:data_to_chance_const} directly defines a lower bound on the number of samples $N$ for a given $p$ and $\alpha$, the lower bound implied by condition \eqref{eq:proof2_N_alpha_p1} in Corollary \ref{thm:data_to_chance_const_p1} can be found numerically for a given $\alpha$.
Because the LHS of \eqref{eq:proof2_N_alpha_p1} decreases exponentially with $N$, it does not require very many samples to satisfy the condition.
For example, $N=36$ is the smallest integer that satisfies condition \eqref{eq:proof2_N_alpha_p1} for the common setting $\alpha=0.01$.

\begin{figure}
     \centering
     \begin{subfigure}[b]{0.6\textwidth}
         \centering
         \includegraphics[width=\textwidth]{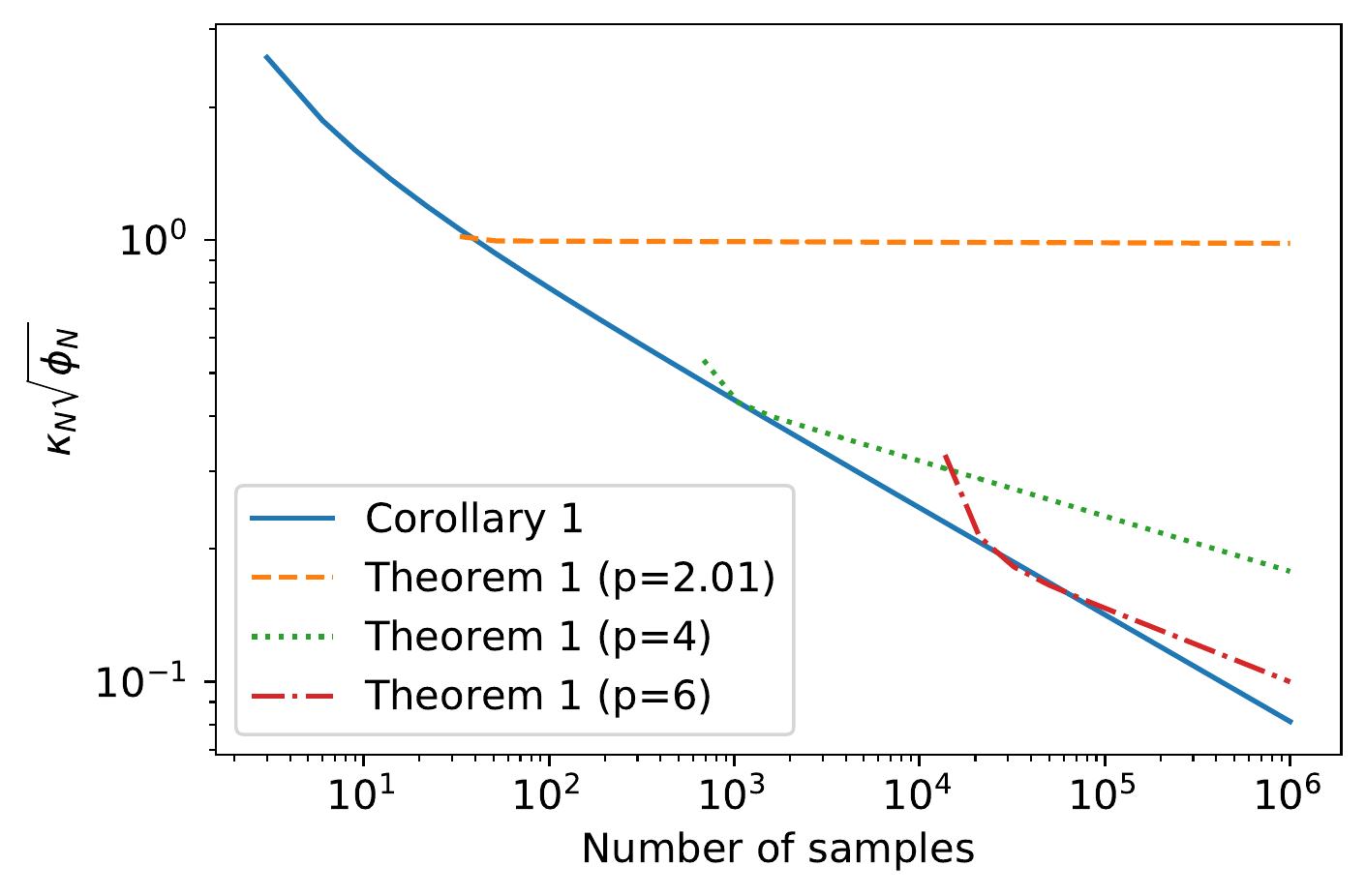}
         \caption{}
         \label{fig:gen_rvs_kappaNsqrtphiN}
     \end{subfigure}
     \hfill
     \begin{subfigure}[b]{0.48\textwidth}
         \centering
         \includegraphics[width=\textwidth]{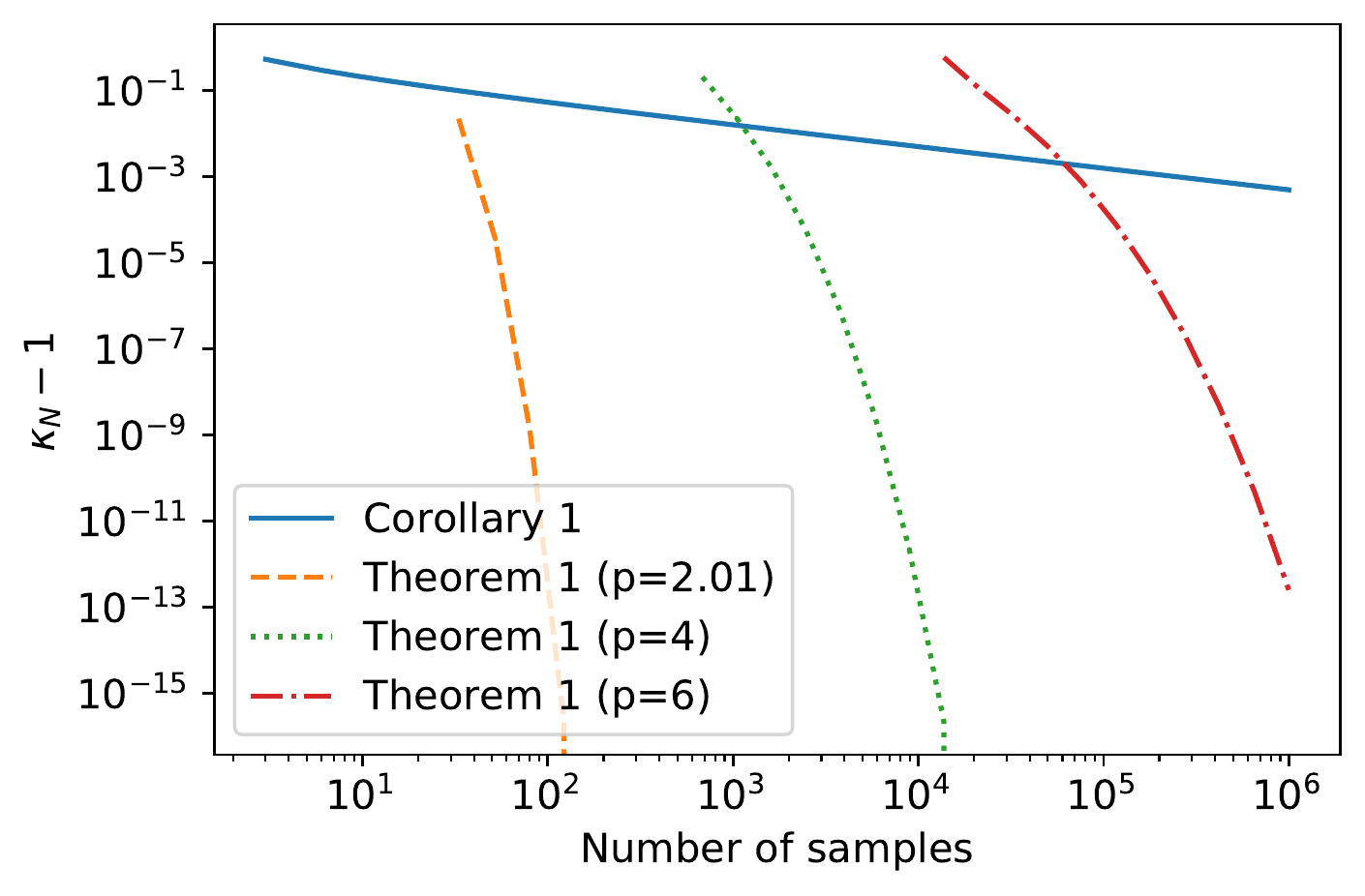}
         \caption{}
         \label{fig:gen_rvs_kappaNminus1}
     \end{subfigure}
     \hfill
     \begin{subfigure}[b]{0.48\textwidth}
         \centering
         \includegraphics[width=\textwidth]{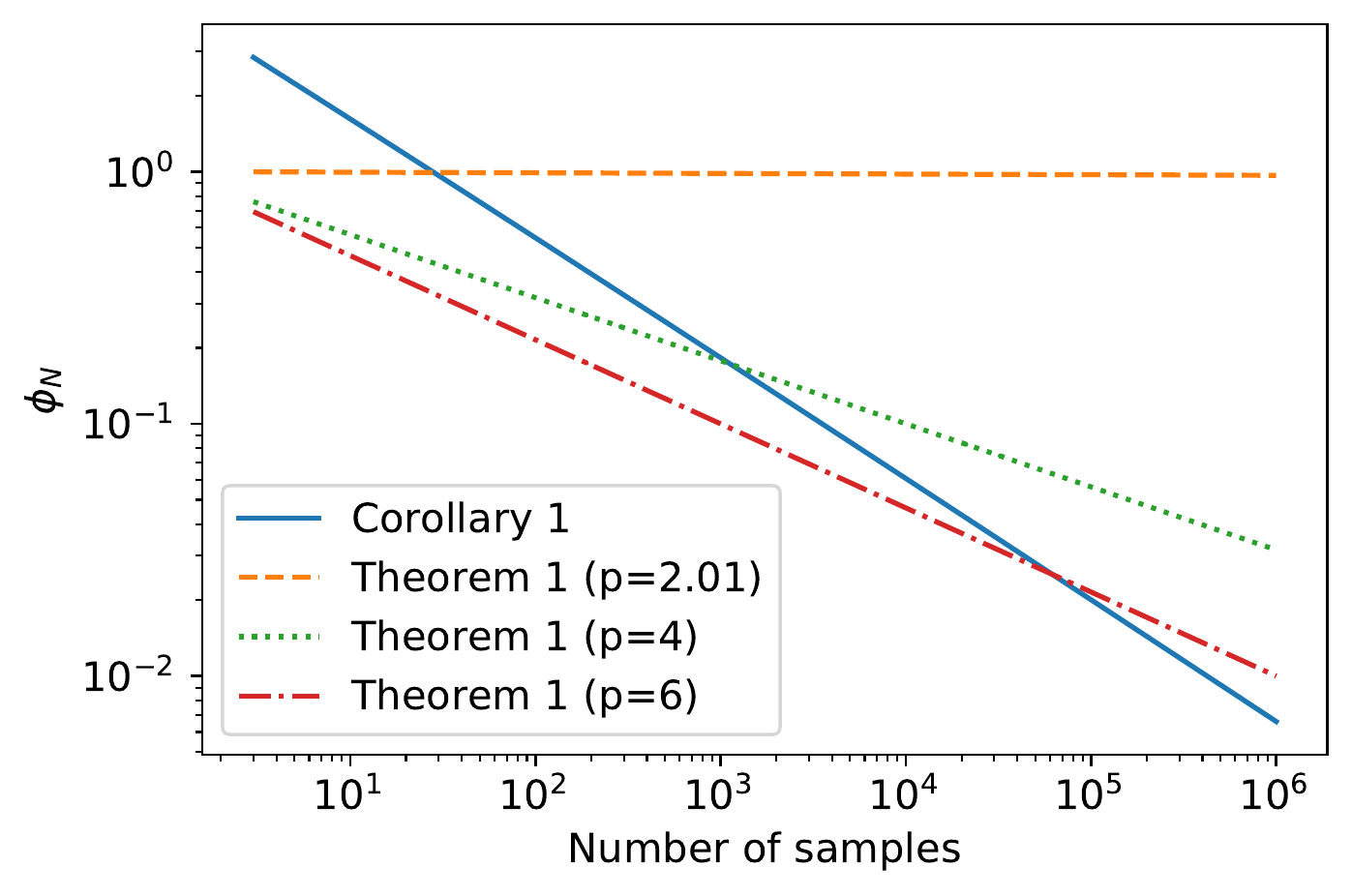}
         \caption{}
         \label{fig:gen_rvs_phiN}
     \end{subfigure}
    \caption{Comparing Corollary \ref{thm:data_to_chance_const_p1} with Theorem \ref{thm:data_to_chance_const} under various settings of $p$ for values of (a) $\kappa_N\sqrt{\phi_N}$, (b) $\kappa_N-1$, and (c) $\phi_N$ vs. the number of samples $N$ when $\alpha=0.1$.}
    \label{fig:gen_rvs}
\end{figure}

The speed at which the deterministic constraint \eqref{eq:prob_cc_distrob_data} approaches the distributionally robust chance constraint \eqref{eq:prob_cc_distrob_true} depends on how quickly $\kappa_N$ approaches 1, $\phi_N$ approaches 0, and $\kappa_N\sqrt{\phi_N}$ approaches 0 with the number of samples $N$ approaching infinity.
By the law of large numbers, we also have that $\hat{\boldsymbol{\mu}}_N$ approaches $\boldsymbol{\mu}^*$, and $\hat{\boldsymbol{\Sigma}}_N$ approaches $\boldsymbol{\Sigma}^*$ as the number of samples $N$ approaches infinity.
When $\kappa_N=1$, $\phi_N=0$, $\hat{\boldsymbol{\mu}}_N=\boldsymbol{\mu}^*$, $\hat{\boldsymbol{\Sigma}}_N=\boldsymbol{\Sigma}^*$, and constraints \eqref{eq:prob_cc_distrob_data_y1} and \eqref{eq:prob_cc_distrob_data_y1} are tight, then \eqref{eq:prob_cc_distrob_data} is equivalent to \eqref{eq:calafiore3p1all}, and thus it is also equivalent to \eqref{eq:prob_cc_distrob_true} due to Lemma \ref{thm:calafiore3p1all}.
In Figure \ref{fig:gen_rvs}, we set $\alpha=0.1$ and compare Corollary \ref{thm:data_to_chance_const_p1} to Theorem \ref{thm:data_to_chance_const} with regards to their effect on \eqref{eq:prob_cc_distrob_data} versus the number of samples.
We can see that although $\kappa_N$ for Theorem \ref{thm:data_to_chance_const} approaches 1 much faster than that for Corollary \ref{thm:data_to_chance_const_p1}, the corresponding $\phi_N$ approaches 0 much slower.
The main advantage to using Corollary \ref{thm:data_to_chance_const_p1} is that for $\kappa_N\sqrt{\phi_N}$ it almost perfectly traces the optimal frontier of Theorem \ref{thm:data_to_chance_const} with a specific $p$ that is given as a function of the number of samples.
Thus, Corollary \ref{thm:data_to_chance_const_p1} strikes an almost optimal balance between the asymptotic behavior of $\kappa_N$ and $\phi_N$.
Later, in Section \ref{sec:betting}, we use simulations to show how using Corollary \ref{thm:data_to_chance_const_p1} can simultaneously give advantageous properties compared to setting $p$ just above 2 and setting a large value of $p$ in Theorem \ref{thm:data_to_chance_const}.

\subsection{Relaxing the Estimated Support Assumption}
\label{sec:relax_support}

There can be many scenarios where it is impractical to use an estimated support that assumes it contains the true support.
This can happen when the true support is very large compared to where most of the samples are located or when the true support is even unbounded.

We can adjust these results to allow for a more relaxed version of the estimated support.
Instead of requiring that the estimated support $\hat{\mathcal{S}}$ contains the true support $\mathcal{S}^*$, let $\hat{\mathcal{S}}_{1-\epsilon}$ be a set that contains a realization of $\mathbf{a}\sim P^*$ with probability at least $1-\epsilon$.
We formally define $\hat{\mathcal{S}}_{1-\epsilon}$ as any set that satisfies the following condition:
\begin{align}
    \text{Pr}(\mathbf{a}\in\hat{\mathcal{S}}_{1-\epsilon} |\mathbf{a}\sim P^*) \geq 1-\epsilon \label{eq:relax_support_cc}
\end{align}
which is 
a probabilistic uncertainty set. 
Note that if $\epsilon=0$, then it is equivalent to the previously defined estimated support $\hat{\mathcal{S}}$.

Our adjustment to the previously described methods is to: (i) replace the estimated support $\hat{\mathcal{S}}$ with the probabilistic uncertainty set $\hat{\mathcal{S}}_{1-\epsilon}$;
(ii) use only a subset of the samples where $\mathbf{a}_i\in \hat{\mathcal{S}}_{1-\epsilon}$ to form the sample mean and covariance; and (iii) replace the scalar $\alpha$ with $\tilde{\alpha}:=\frac{\alpha-\epsilon}{1-\epsilon}$, where it is assumed that $\epsilon < \alpha$.
Under this scenario, it is assumed that the practitioner knows beforehand a value of $\epsilon$ for the estimated support $\hat{\mathcal{S}}_{1-\epsilon}$ that satisfies condition \eqref{eq:relax_support_cc}.
The following proposition formalizes this adjustment to the original method so that the results of Theorem \ref{thm:data_to_chance_const} apply.
\begin{proposition}
    \label{thm:Ptrunc_thm1}
    Assume that the estimated support $\hat{\mathcal{S}}_{1-\epsilon}$ satisfies \eqref{eq:relax_support_cc} for some $\epsilon<\alpha\in(0,1)$ and that at least $N$ of the independently drawn samples are contained in $\hat{\mathcal{S}}_{1-\epsilon}$, i.e., $\mathbf{a}_i\in\hat{\mathcal{S}}_{1-\epsilon}:\forall i\in\{1,\dots,N\}$.
    Let there exist a scalar $p > 2$ such that:
        \begin{align}
            N & >\left(2+\sqrt{2\ln(4/\tilde{\alpha})}\right)^p 
        \end{align}
    where $\tilde{\alpha}:=\frac{\alpha-\epsilon}{1-\epsilon}$.
    Define the following scalars that depend on $N$, $\tilde{\alpha}$, and $p$:
    \begin{subequations}
        \begin{align}
            \kappa_N & := \left(1-\frac{4}{\tilde{\alpha}}\exp\left(-(N^\frac{1}{p}-2)^2/2\right)\right)^{-\frac{1}{2}} \\
            \phi_N & := N^{\left(\frac{1}{p}-\frac{1}{2}\right)}.
        \end{align}
    \end{subequations}
    Then, the convex deterministic constraints \eqref{eq:prob_cc_distrob_data}
    represent a conservative approximation of the distributionally robust chance constraint \eqref{eq:prob_cc_distrob_true}. In addition,  the constraints \eqref{eq:prob_cc_distrob_data} asymptotically approach \eqref{eq:prob_cc_distrob_true} as $N\rightarrow\infty$.
\end{proposition}

\proof
    By applying Theorem \ref{thm:data_to_chance_const} with $\tilde{\alpha}$ for the probability distribution $P^*$ truncated by $\hat{\mathcal{S}}_{1-\epsilon}$, we have that:
    \begin{align}
        \inf_{P\in\mathcal{P}(\boldsymbol{\mu}^*(\hat{\mathcal{S}}_{1-\epsilon}),\boldsymbol{\Sigma}^*(\hat{\mathcal{S}}_{1-\epsilon}),\hat{\mathcal{S}}_{1-\epsilon})}\text{Pr}\left(\mathbf{a}^\intercal\mathbf{x}\leq 0|\mathbf{a}\sim P \right) \geq 1- \tilde{\alpha}
    \end{align}
    where $\boldsymbol{\mu}^*(\hat{\mathcal{S}}_{1-\epsilon})$ is the mean, and $\boldsymbol{\Sigma}^*(\hat{\mathcal{S}}_{1-\epsilon})$ is the covariance of the truncated distribution, and $\mathcal{P}(\cdot,\cdot,\hat{\mathcal{S}}_{1-\epsilon})$ is the family of probability distributions truncated by $\hat{\mathcal{S}}_{1-\epsilon}$ for a given mean and covariance.
    
    The rest of the proof is similar to that of Lemma \ref{thm:confregion} and is a simple application of the law of total probability to show that this distributionally robust chance constraint implies the distributionally robust chance constraint for the probability distribution $P^*$ with $\alpha$.
    \begin{align}
        \inf_{P\in\mathcal{P}(\boldsymbol{\mu}^*,\boldsymbol{\Sigma}^*)}\text{Pr}\left(\mathbf{a}^\intercal\mathbf{x}\leq 0|\mathbf{a}\sim P \right) & = \inf_{P\in\mathcal{P}(\boldsymbol{\mu}^*,\boldsymbol{\Sigma}^*)}\{\text{Pr}(\mathbf{a}^\intercal\mathbf{x}\leq 0|\mathbf{a}\sim P,\mathbf{a}\in \hat{\mathcal{S}}_{1-\epsilon})\text{Pr}( \mathbf{a}\in \hat{\mathcal{S}}_{1-\epsilon}|\mathbf{a}\sim P ) \nonumber \\
        & \quad\quad\quad\quad\quad\quad\quad\quad + \text{Pr}(\mathbf{a}^\intercal\mathbf{x}\leq 0|\mathbf{a}\sim P, \mathbf{a}\notin \mathcal{S}_{1-\epsilon} ) \text{Pr}( \mathbf{a}\notin \hat{\mathcal{S}}_{1-\epsilon}|\mathbf{a}\sim P )\} \nonumber \\
        & \geq \inf_{P\in\mathcal{P}(\boldsymbol{\mu}^*,\boldsymbol{\Sigma}^*)}\{\text{Pr}(\mathbf{a}^\intercal\mathbf{x}\leq 0|\mathbf{a}\sim P,\mathbf{a}\in \hat{\mathcal{S}}_{1-\epsilon} )\text{Pr}( \mathbf{a}\in \hat{\mathcal{S}}_{1-\epsilon}|\mathbf{a}\sim P )\} \nonumber \\
        & = \inf_{P\in\mathcal{P}(\boldsymbol{\mu}^*,\boldsymbol{\Sigma}^*)}\{\text{Pr}(\mathbf{a}^\intercal\mathbf{x}\leq 0|\mathbf{a}\sim P,\mathbf{a}\in \hat{\mathcal{S}}_{1-\epsilon} )\}(1-\epsilon) \nonumber \\
        & = \inf_{P\in\mathcal{P}(\boldsymbol{\mu}^*(\hat{\mathcal{S}}_{1-\epsilon}),\boldsymbol{\Sigma}^*(\hat{\mathcal{S}}_{1-\epsilon}),\hat{\mathcal{S}}_{1-\epsilon})}\{\text{Pr}\left(\mathbf{a}^\intercal\mathbf{x}\leq 0|\mathbf{a}\sim P \right)\}(1-\epsilon) \nonumber \\
        & \geq (1- \tilde{\alpha})(1-\epsilon) \nonumber \\
        & = \left(1- \frac{\alpha-\epsilon}{1-\epsilon}\right)(1-\epsilon) \nonumber \\
        & = 1-\alpha \nonumber.
    \end{align}
\endproof

Similarly, Corollary \ref{thm:data_to_chance_const_p1} can be adapted for this relaxed case from proof structure of Proposition \ref{thm:Ptrunc_thm1}.

\section{Guarantees for Independent Random Variables}
\label{sec:indep_data-driven_guar}
When the random variables in the vector $\mathbf{a}$ are independent of each other, this information can be used to reduce the conservativeness of the proposed distributionally robust chance-constrained approach. In general, the independence of the random variables makes worst-case realizations of $\mathbf{a}^\intercal\mathbf{x}$ less probable.
The independence assumption also results in decoupling the support and the estimation support along the dimensions, i.e., the support becomes independent intervals.
Let $(\underline{\mathbf{s}}^*,\overline{\mathbf{s}}^*)$ define the true support of $\mathbf{a}$ and $(\underline{\mathbf{s}},\overline{\mathbf{s}})$ define the estimated support such that $\underline{\mathbf{s}}\leq\underline{\mathbf{s}}^*$ and $\overline{\mathbf{s}}\geq\overline{\mathbf{s}}^*$.
Thus, the radius of $\hat{\mathcal{S}}$ can be written as a function of the decision $\mathbf{x}$ as defined by  \eqref{eq:rx_indint}.
Note that the relaxation results from Proposition \ref{thm:Ptrunc_thm1} can also be applied in this case.

Another useful property of independence is that the data sampled for each random variable do not need to be sampled simultaneously or even have the same number of samples.
This property is crucial when the data are collected from independent agents in a distributed system because this means that the random variable data can be measured asynchronously.
Additionally, each agent can preserve some privacy, hence it is only required to report its sample mean $\hat{\mu}_i:=\frac{1}{N}\sum_{j=1}^Na_{j,i}$, its sample variance $\hat{\sigma}^2_i:=\frac{1}{N}\sum_{j=1}^N(a_{j,i}-\hat{\mu}_i)^2$, and the number of measured samples if it is different between agents. 
For simplicity, we assume that there are exactly $N$ samples taken for each random variable, but in practice $N$ can be thought of as the smallest number of samples taken among the independent random variables.

The property of independent random variables allows us to devise two separate distributionally robust methods.
The first method only requires the sample means and the estimated support intervals, whereas the second method additionally requires the use of the sample variances, like in Section \ref{sec:data-driven_gaur}.
Thus, the situation on whether the true probability distribution has large support intervals or large variances will determine which method is closer to the original chance constraint of \eqref{eq:prob_cc_original}.

\subsection{Using Sample Means and Estimated Independent Intervals}

Let $\mathcal{P}_\text{ind}(\boldsymbol{\mu},\underline{\mathbf{s}},\overline{\mathbf{s}})$ be defined as the family of element-wise independent probability distributions with mean $\boldsymbol{\mu}$, which have support bounds $\underline{\mathbf{s}}$ and $\overline{\mathbf{s}}$.
Notice that compared to the family $\mathcal{P}(\boldsymbol{\mu},\boldsymbol{\Sigma})$, we are replacing information from the covariance with the support.  
Thus, the distributionally robust chance constraint defined by the family $\mathcal{P}_\text{ind}(\boldsymbol{\mu}^*,\underline{\mathbf{s}}^*,\overline{\mathbf{s}}^*)$ is:
\begin{align}
    \inf_{P\in\mathcal{P}_\text{ind}(\boldsymbol{\mu}^*,\underline{\mathbf{s}}^*,\overline{\mathbf{s}}^*)}\text{Pr}\left(\mathbf{a}^\intercal\mathbf{x}\leq 0|\mathbf{a}\sim P\right) \geq 1- \alpha  \label{eq:prob_cc_distrob_true_ind}
\end{align}
which conservatively guarantees that the original chance constraint \eqref{eq:prob_cc_original} is satisfied.
If the true mean $\boldsymbol{\mu}^*$ and support $(\underline{\mathbf{s}}^*,\overline{\mathbf{s}}^*)$ are known, the following lemma from \cite{calafiore2006distributionally} gives a convex deterministic constraint that conservatively approximates \eqref{eq:prob_cc_distrob_true_ind}.

\begin{lemma}
    \label{thm:calafiore3p2}
    \cite{calafiore2006distributionally}(Lemma 3.2) The distributionally robust chance constraint \eqref{eq:prob_cc_distrob_true_ind} for $\alpha\in(0,1)$ is satisfied if:
    \begin{align}\label{eq:calafiore3p2_constraint}
        \boldsymbol{\mu}^{*\intercal}\mathbf{x}  + \sqrt{\frac{1}{2}\ln\left(\frac{1}{\alpha}\right)}\|\mathbf{S}^*\mathbf{x}\|_2 \leq 0
    \end{align}
    where $\mathbf{S}^*:=\text{diag}(\overline{\mathbf{s}}^*-\underline{\mathbf{s}}^*)$.
\end{lemma}

\begin{figure}
    \centering
    \includegraphics[width=0.5\textwidth]{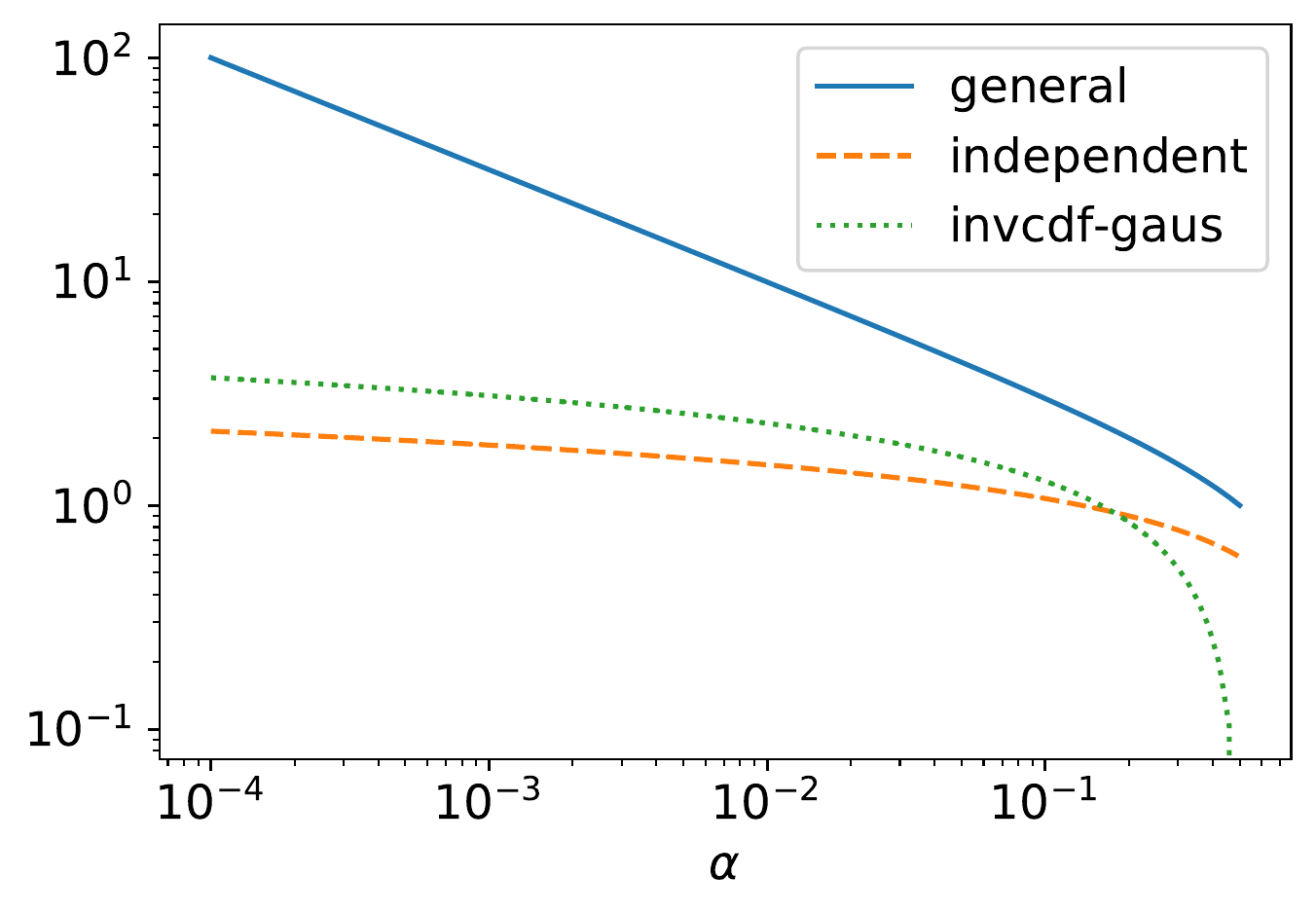}
    \caption{Deterministic chance constraint equivalents and approximations vs. $\alpha$ for $\sqrt{\frac{1-\alpha}{\alpha}}$ (``general"), $\sqrt{\frac{1}{2}\ln\left(\frac{1}{\alpha}\right)}$ (``independent") and the inverse cumulative distribution function of the Gaussian distribution $\Phi^{-1}(1-\alpha)$.  If $\sqrt{\frac{1-\alpha}{\alpha}}$ were replaced by $\Phi^{-1}(1-\alpha)$ in Lemma \ref{thm:calafiore3p1all}, then the lemma would instead give the deterministic equivalent of the Gaussian chance constraint (See \cite{calafiore2006distributionally}, Section 2.1)}
    \label{fig:alpha_v_const}
\end{figure}

Note that compared to Lemma \ref{thm:calafiore3p1all}, Lemma \ref{thm:calafiore3p2} uses a diagonal matrix of the support bounds instead of the covaraince matrix. Note that although Lemma \ref{thm:calafiore3p1all} provided an equivalent reformulation, this lemma provides only a conservative approximation of its associated distributionally robust chance constraint.
On the other hand, the value of $\sqrt{\frac{1-\alpha}{\alpha}}$ from Lemma 
\ref{thm:calafiore3p1all} grows faster as $\alpha$ goes toward zero than its counterpart $\sqrt{\frac{1}{2}\ln\left(\frac{1}{\alpha}\right)}$ in Lemma \ref{thm:calafiore3p2}.
Figure \ref{fig:alpha_v_const} compares these functions versus $\alpha$, which shows the numerical advantage of assuming independence when possible.

Similar to the case without the independence assumption, we form a convex data-driven deterministic surrogate constraint to replace \eqref{eq:prob_cc_distrob_true_ind} in an optimization problem such as \eqref{eq:prob_cc_distrob_data} using the sample mean $\hat{\boldsymbol{\mu}}_N$ and the estimated support $\mathbf{S}:=\text{diag}(\overline{\mathbf{s}}-\underline{\mathbf{s}})$:
\begin{align}
    \hat{\boldsymbol{\mu}}_N^\intercal\mathbf{x} + \left(\frac{1}{2}\phi_N + \sqrt{\frac{1}{2}\ln\left(\frac{1}{\alpha}\right)+\nu_N}\right)\|\mathbf{S}\mathbf{x}\|_1 & \leq 0. \label{eq:prob_cc_distrob_data_ind}
\end{align}
The following theorem proves that with specific scalars $(\phi_N,\nu_N)$, and under a condition on the number of samples $N$, the deterministic constraints \eqref{eq:prob_cc_distrob_data_ind} are indeed a conservative approximation of \eqref{eq:prob_cc_distrob_true_ind}.

\begin{proposition}
    \label{thm:data_to_chance_const_ind}
    Assume that the estimated element-wise independent support intervals $(\underline{\mathbf{s}},\overline{\mathbf{s}})$ contain the true support intervals $(\underline{\mathbf{s}}^*,\overline{\mathbf{s}}^*)$, i.e., $\underline{\mathbf{s}}\leq\underline{\mathbf{s}}^*$ and $\overline{\mathbf{s}}\geq\overline{\mathbf{s}}^*$.
    Let there exist a scalar $p$ such that:
    \begin{subequations}
        \label{eq:proof2_DRCC_conditions_ind}
        \begin{align}
            N & >\left(2+\sqrt{2\ln(1/\alpha)}\right)^p \label{eq:proof2_DRCC_conditions_ind_a} \\
            p & > 0
        \end{align}
    \end{subequations}
    where $\alpha\in(0,1)$.
    Define the following scalars that depend on $N$, $\alpha$, and $p$:
    \begin{subequations}
        \label{eq:proof2_DRCC_constants_ind}
        \begin{align}
            \nu_N & := \frac{1}{2}\ln{\left(1 + \frac{1-\alpha}{\alpha\exp{\left((N^{\frac{1}{p}}-2)^2/2\right)}-1}\right)} \\
            \phi_N & := N^{\left(\frac{1}{p}-\frac{1}{2}\right)}.
        \end{align}
    \end{subequations}
    Then, the convex deterministic constraint \eqref{eq:prob_cc_distrob_data_ind}
    is a conservative approximation of the  distributionally robust chance constraint \eqref{eq:prob_cc_distrob_true_ind}.
\end{proposition}

The proof is structured very similarly to that of Theorem \ref{thm:data_to_chance_const} and is in Appendix \ref{prf:data_to_chance_const_ind}.
Instead of Lemma \ref{thm:calafiore4p1}, it uses a simpler version from \cite{shawe2003estimating}, which considers only the mean.

\begin{lemma}
    \label{thm:shawe}\cite{shawe2003estimating}(Theorem 3) Let $\mathbf{v}_1,\dots\mathbf{v}_N$ be $N$ independent samples drawn from the probability distribution $P^*$ with support $\mathcal{S}^*$ for $\delta\in(0,1)$.
    Define $r^*$ as the radius of $\mathcal{S}^*$, i.e., $r^*:=\frac{1}{2}\sup_{\{\mathbf{v}_1,\mathbf{v}_2\}\in\mathcal{S}^*}\|\mathbf{v}_1-\mathbf{v}_2\|_2$.
    Then:
    \begin{align}
        \|\boldsymbol{\mu}^*-\hat{\boldsymbol{\mu}}_N\|_2 & \leq \frac{r^*}{\sqrt{N}}\left(2+\sqrt{2\ln(1/\delta)}\right)
    \end{align}
    with probability at least $1-\delta$.
\end{lemma}

Similar to Corollary \ref{thm:data_to_chance_const_p1}, we simplify the form of $(\phi_N,\nu_N)$ by setting $p$ to be a function of the number of samples $N$ and $\alpha$ as follows.
\begin{corollary}
    \label{thm:data_to_chance_const_p1_ind}
    Theorem \ref{thm:data_to_chance_const_ind} remains true if the conditions in \eqref{eq:proof2_DRCC_conditions_ind} are replaced with:
    \begin{align}
        N \geq 2 \label{eq:proof2_N_alpha_p1_ind}
    \end{align}
    and the scalar definitions in \eqref{eq:proof2_DRCC_constants_ind} are replaced with:
    \begin{subequations}\label{eq:proof2_DRCC_constants_p1_ind}
        \begin{align}
            \nu_N & := \frac{1}{2}\ln{\left(1 + \frac{1-\alpha}{\sqrt{N}-1}\right)} \\
            \phi_N & := \frac{2+\sqrt{2\ln (\sqrt{N}/\alpha)}}{\sqrt{N}}.
        \end{align}
    \end{subequations}
\end{corollary}
\proof
    Let $p:=\log_{2+\sqrt{2\ln (\sqrt{N}/\alpha)}}(N)>0$, then we can write $N^{\frac{1}{p}} = 2+\sqrt{2\ln (\sqrt{N}/\alpha)}$. By simple substitution in \eqref{eq:proof2_DRCC_constants_ind}, we obtain the expression in \eqref{eq:proof2_DRCC_constants_p1_ind}. In addition, this specific setting of $p$ ensures that $N^{\frac{1}{p}} \geq 2+\sqrt{2\ln (1/\alpha)}$ for $N\geq 2$, which completes the conditions from Theorem \ref{thm:data_to_chance_const_ind} and thus proves the corollary.
\endproof

\begin{figure}
     \centering
     \begin{subfigure}[b]{0.48\textwidth}
         \centering
         \includegraphics[width=\textwidth]{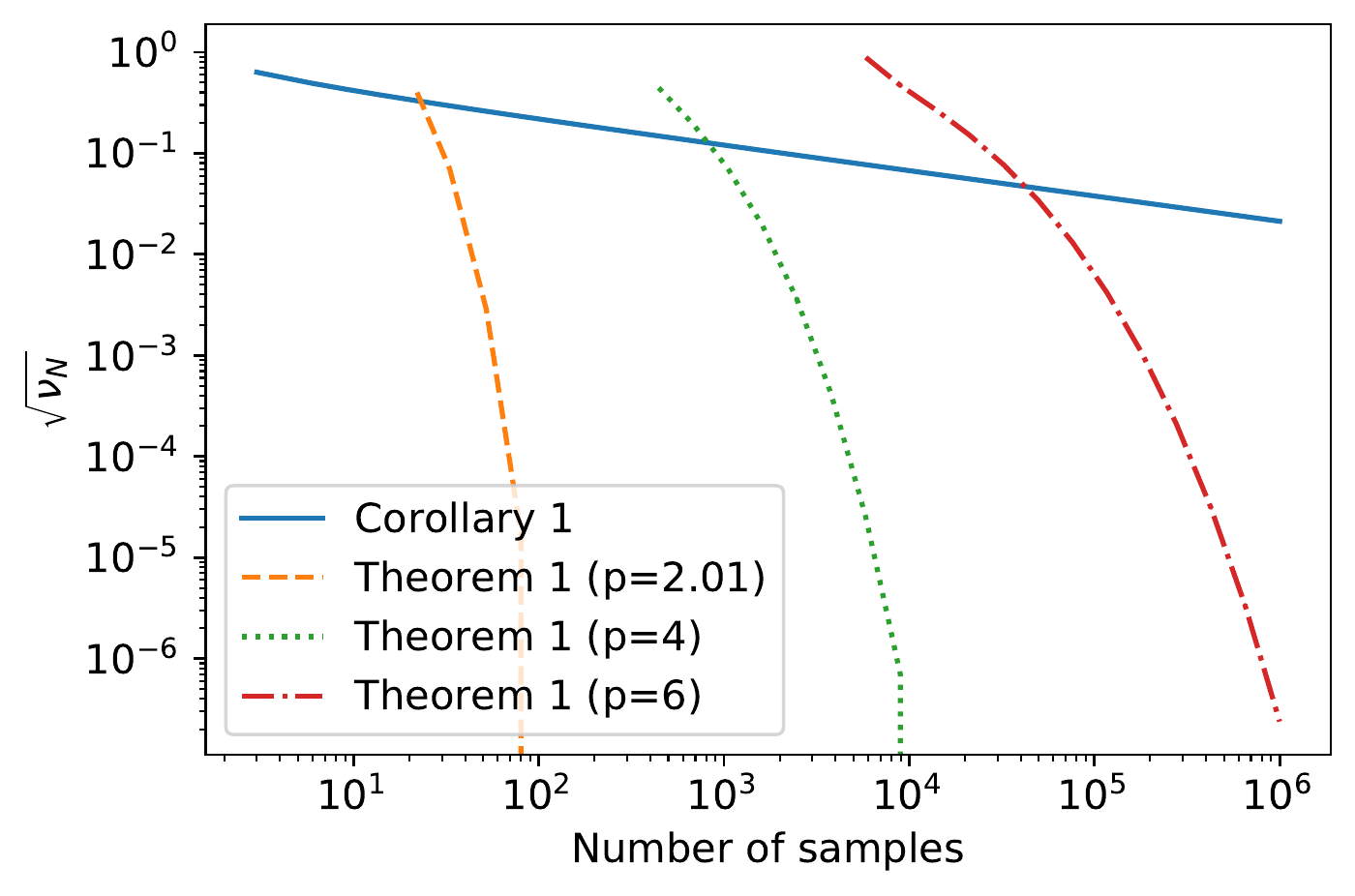}
         \caption{}
         \label{fig:ind_rvs_sqrtnuN}
     \end{subfigure}
     \hfill
     \begin{subfigure}[b]{0.48\textwidth}
         \centering
         \includegraphics[width=\textwidth]{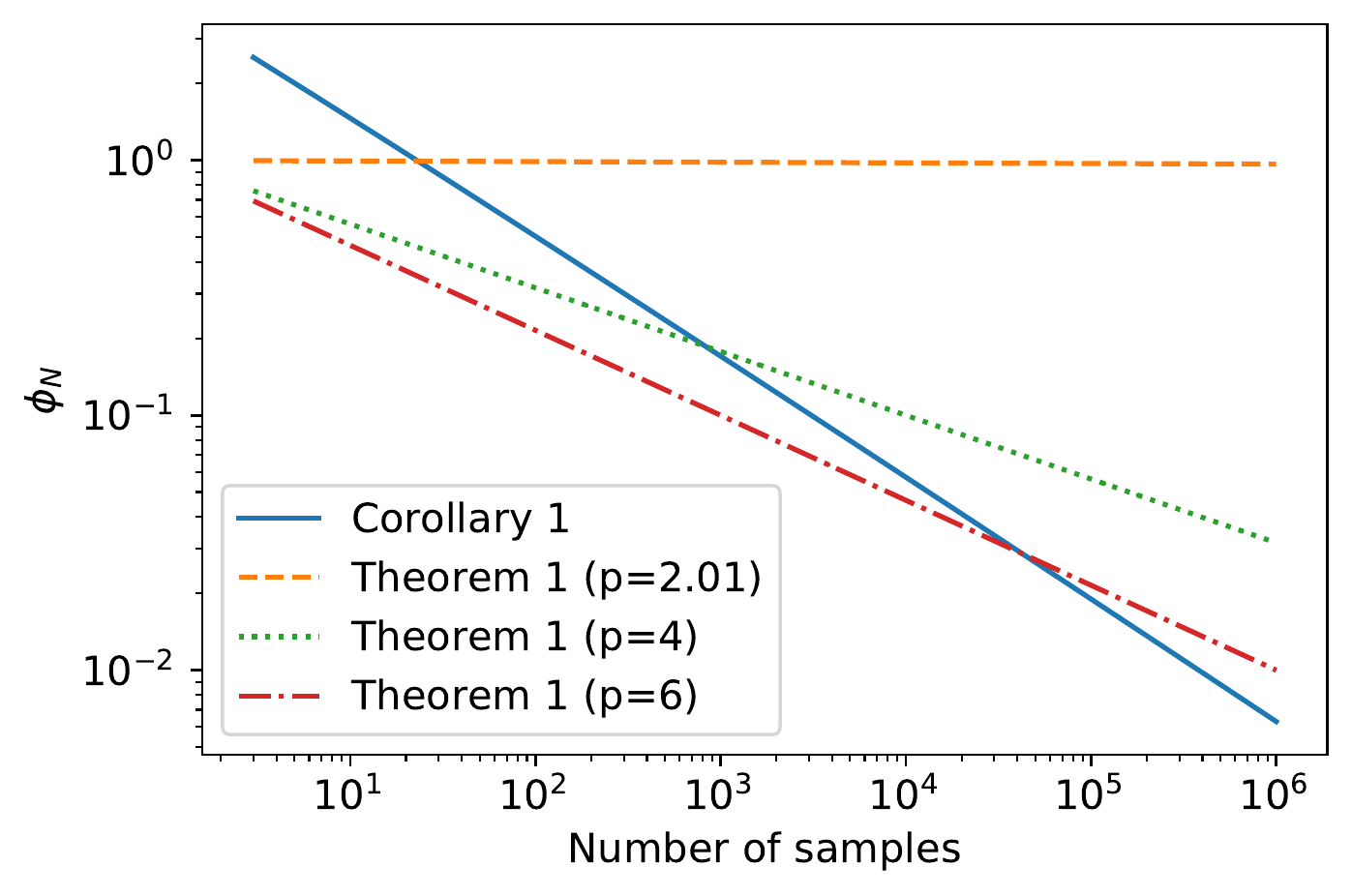}
         \caption{}
         \label{fig:ind_rvs_phiN}
     \end{subfigure}
    \caption{Comparing Corollary \ref{thm:data_to_chance_const_p1} with Theorem \ref{thm:data_to_chance_const} under various settings of $p$ for values of (a) $\kappa_N\sqrt{\phi_N}$, (b) $\kappa_N-1$, and (c) $\phi_N$ vs. the number of samples $N$ when $\alpha=0.1$.}
    \label{fig:ind_rvs}
\end{figure}

The speed at which the deterministic constraint \eqref{eq:prob_cc_distrob_data_ind} approaches Lemma \ref{thm:calafiore3p2} depends on how quickly $\sqrt{\nu_N}$ and $\phi_N$ approach 0.
In Figure \ref{fig:ind_rvs}, we set $\alpha=0.1$ and compare Corollary \ref{thm:data_to_chance_const_p1_ind} to Proposition \ref{thm:data_to_chance_const_ind} with regards to their effect on \eqref{eq:prob_cc_distrob_data_ind} versus the number of samples.
We can see that although $\sqrt{\nu_N}$ for Proposition \ref{thm:data_to_chance_const_ind} approaches 0 much faster than that for Corollary \ref{thm:data_to_chance_const_p1_ind}, its $\phi_N$ approaches 0 much slower.
The main advantage to using Corollary \ref{thm:data_to_chance_const_p1_ind} is that $\sqrt{\nu_N}$ and $\phi_N$ approach 0 at approximately the same rate, which makes an almost optimal balance between the asymptotic behavior of $\kappa_N$ and $\phi_N$ alleviating the requirement to tune $p$.

\subsection{Using Sample Means, Sample Variances, and Estimated Independent Intervals}

With a slight abuse of notation, we say that $\mathcal{P}_\text{ind}(\boldsymbol{\mu},\boldsymbol{\sigma}^2)$ is defined as the family of element-wise independent probability distributions with mean $\boldsymbol{\mu}$ and variance $\boldsymbol{\sigma}^2$.
Thus, the distributionally robust chance constraint defined by the family $\mathcal{P}_\text{ind}(\boldsymbol{\mu}^*,\boldsymbol{\sigma}^{2*})$ is:
\begin{align}
    \inf_{P\in\mathcal{P}_\text{ind}(\boldsymbol{\mu}^*,\boldsymbol{\sigma}^{2*})}\text{Pr}\left(\mathbf{a}^\intercal\mathbf{x}\leq 0|\mathbf{a}\sim P\right) \geq 1- \alpha  \label{eq:prob_cc_distrob_true_indvar}
\end{align}
which conservatively guarantees that the original chance constraint \eqref{eq:prob_cc_original} is satisfied.
If the true mean $\boldsymbol{\mu}^*$ and variance $\boldsymbol{\sigma}^{2*}$ are known, the following lemma gives a convex deterministic constraint that conservatively approximates \eqref{eq:prob_cc_distrob_true_indvar}, which is just a special case of Lemma \ref{thm:calafiore3p1all}.

\begin{lemma}
    \label{thm:calafiore3p1all_indvar}
    The distributionally robust chance constraint \eqref{eq:prob_cc_distrob_true_indvar} for $\alpha\in(0,1)$ is equivalent to:
    \begin{align}\label{eq:prob_cc_distrob_true_indvar_equiv}
        \boldsymbol{\mu}^{*\intercal}\mathbf{x} + \sqrt{\frac{1-\alpha}{\alpha}}\|\mathbf{D}^*\mathbf{x}\|_2 \leq 0
    \end{align}
    where $\mathbf{D}^*$ is a diagonal matrix such that $D^*_{ii}=\sqrt{\sigma^{2*}_i}:\forall i\in\{1,\dots,n\}$.
\end{lemma}
\proof
    Assuming that the underlying distribution of $\mathbf{a}$ is element-wise independent, i.e., all the off-diagonal elements of $\boldsymbol{\Sigma}$ are zero, we get \eqref{eq:prob_cc_distrob_true_indvar_equiv} as a special case of Lemma \ref{thm:calafiore3p1all}~\cite{calafiore2006distributionally}. 
\endproof

With the sample mean $\hat{\boldsymbol{\mu}}_N$, sample variance $\hat{\boldsymbol{\sigma}}^2_N$, and the radius $\frac{1}{2}\|\mathbf{S}\mathbf{x}\|_1$ of the estimated support $(\underline{\mathbf{s}},\overline{\mathbf{s}})$, we can replace the distributionally robust chance constraint \eqref{eq:prob_cc_distrob_true_indvar} in an optimization problem with the following deterministic convex constraint:
\begin{subequations}
    \label{eq:prob_cc_distrob_data_indvar}
    \begin{align}
        \hat{\boldsymbol{\mu}}_N^\intercal\mathbf{x} + \frac{1}{2}\phi_N \|\mathbf{S}\mathbf{x}\|_1 + \kappa_N\sqrt{\frac{1-\alpha}{\alpha}}\|\mathbf{y}\|_2 & \leq 0 \\
        \|\hat{\mathbf{D}}_N\mathbf{x}\|_2 & \leq y_1 \\
        \sqrt{\frac{1}{2}\phi_N}\|\mathbf{S}\mathbf{x}\|_1 & \leq y_2.
    \end{align}
\end{subequations}
where $\mathbf{y}\in\mathbb{R}^2$ is a vector of auxiliary variables, $\hat{\mathbf{D}}_N$ is a diagonal matrix that collects the sample variances such that $\hat{D}_{N,ii}=\sqrt{\hat{\sigma}^{2}_i}:\forall i\in\{1,\dots,n\}$, and $(\phi_N,\kappa_N)$ are scalars that depend on the number of samples $N$ and $\alpha$.
This deterministic constraint is a conservative approximation that guarantees \eqref{eq:prob_cc_distrob_true_indvar} will be satisfied. In addition, it asymptotically approaches \eqref{eq:prob_cc_distrob_true_indvar} as the number of samples $N$ approaches infinity.
The following theorem is identical to Theorem \ref{thm:data_to_chance_const}, but it specifically considers the case of independent random variables.

\begin{proposition}
    \label{thm:data_to_chance_const_indvar}
    Assume that the estimated element-wise independent support intervals $(\underline{\mathbf{s}},\overline{\mathbf{s}})$ contain the true support intervals $(\underline{\mathbf{s}}^*,\overline{\mathbf{s}}^*)$, i.e., $\underline{\mathbf{s}}\leq\underline{\mathbf{s}}^*$ and $\overline{\mathbf{s}}\geq\overline{\mathbf{s}}^*$.
    Let there exist a scalar $p$ such that:
    \begin{subequations}
        \label{eq:proof2_DRCC_conditions_indvar}
        \begin{align}
            N & >\left(2+\sqrt{2\ln(4/\alpha)}\right)^p \label{eq:proof2_DRCC_conditions_a_indvar} \\
            p & > 2
        \end{align}
    \end{subequations}
    where $\alpha\in(0,1)$.
    Define the following scalars that depend on $N$, $\alpha$, and $p$:
    \begin{subequations}
        \label{eq:proof2_DRCC_constants_indvar}
        \begin{align}
            \kappa_N & := \left(1-\frac{4}{\alpha}\exp\left(-(N^\frac{1}{p}-2)^2/2\right)\right)^{-\frac{1}{2}} \\
            \phi_N & := N^{\left(\frac{1}{p}-\frac{1}{2}\right)}.
        \end{align}
    \end{subequations}
    Then, the convex deterministic constraint \eqref{eq:prob_cc_distrob_data_indvar}
    is a conservative approximation of the following distributionally robust chance constraint \eqref{eq:prob_cc_distrob_true_indvar} and asymptotically approaches it as $N\rightarrow\infty$.
\end{proposition}

The proof is exactly the same as that of Theorem \ref{thm:data_to_chance_const} except that Lemma \ref{thm:calafiore3p1all_indvar} replaces Lemma \ref{thm:calafiore3p1all}.
It is necessary to explicitly state this proposition becasue of the technicality that Theorem \ref{thm:data_to_chance_const} assumes that sample covariances between elements of $\mathbf{a}$ are calculated from the data, whereas here we assume that they are all zero by only calculating the individual variance for each element of $\mathbf{a}$.

Likewise, we also have the following corollary.
\begin{corollary}
    \label{thm:data_to_chance_const_p1_indvar}
    Theorem \ref{thm:data_to_chance_const_indvar} remains true if the conditions in \eqref{eq:proof2_DRCC_conditions_indvar} are replaced with:
    \begin{align}
        \sqrt{\frac{16N}{\exp{\left(\left(\sqrt{N}-2\right)^2\right)}}} < \alpha \label{eq:proof2_N_alpha_p1_indvar}
    \end{align}
    and the scalar definitions in \eqref{eq:proof2_DRCC_constants_indvar} are replaced with:
    \begin{subequations}
        \begin{align}
            \kappa_N & := \sqrt{\frac{\sqrt{N}}{\sqrt{N}-1}} \\
            \phi_N & := \frac{2+\sqrt{2\ln (4\sqrt{N}/\alpha)}}{\sqrt{N}}.
        \end{align}
    \end{subequations}
\end{corollary}

Again, the proof follows the same argument as Corollary \ref{thm:data_to_chance_const_p1}.

\begin{remark}
    With the use of the sample mean $\hat{\boldsymbol{\mu}}_N$, sample variance $\hat{\boldsymbol{\sigma}}^2_N$, and the estimated support $(\underline{\mathbf{s}},\overline{\mathbf{s}})$, both deterministic constraints \eqref{eq:prob_cc_distrob_data_ind} and \eqref{eq:prob_cc_distrob_data_indvar} could be used simultaneously to conservatively surrogate
    the following distributionally robust chance constraint:
    \begin{align}
        \inf_{P\in\mathcal{P}_\text{ind}(\boldsymbol{\mu}^*,\boldsymbol{\sigma}^{2*},\underline{\mathbf{s}}^*,\overline{\mathbf{s}}^*)}\text{Pr}\left(\mathbf{a}^\intercal\mathbf{x}\leq 0|\mathbf{a}\sim P\right) \geq 1- \alpha.
    \end{align}
    Constraint \eqref{eq:prob_cc_distrob_data_ind} requires a minimum of two samples, whereas \eqref{eq:prob_cc_distrob_data_indvar} requires a number of samples that are considerably larger depending on $\alpha$.
    On the other hand, as the number of samples $N$ increases, \eqref{eq:prob_cc_distrob_data_indvar} converges to the constraint \eqref{eq:prob_cc_distrob_true_indvar_equiv} as if the mean and variance were known a priori; however, the data-driven constraint \eqref{eq:prob_cc_distrob_data_ind} approaches the constraint \eqref{eq:calafiore3p2_constraint} only as close as the estimated support $(\underline{\mathbf{s}},\overline{\mathbf{s}})$ is to the true support $(\underline{\mathbf{s}}^*,\overline{\mathbf{s}}^*)$.
    Therefore, both data-driven constraints should be used separately in an optimization problem when possible. Then, the solution that gives the least conservative objective function value can be selected.
\end{remark}

\section{Numerical Evaluation: Allocating Bets among Correlated Wagers}
\label{sec:betting}
In this section, we evaluate the proposed distributionally robust chance constraint method using a simple application example. This application showcases a problematic probability distribution in that all its realizations are located at the boundary of its support.
Because solving the optimization problem has the largest time complexity and our focus is on large sample sizes, we compare our approach against the data-driven methods found in the literature that do not increase the size of the optimization problem with the number of samples (see top half of Table \ref{tab:cc_complexity}).

\subsection{Setup}
Suppose that a sports better places a bet on the outcome of a game but has several available options of how to bet on the game based on how many points a particular team needs to win by, called a point spread;
however, the outcomes of the available betting options for the same game are strongly correlated.
For example, one wager could be that Team A wins the game by at least 1 point, whereas another wager is that Team A wins by at least 2 points.
Clearly, the outcomes of these options are dependent on each other because the first wager will always be true if the second one is true.

Let $x_i$ be the fraction of the bankroll placed on wager $i$ and $a_i\in\{-1,\overline{a}_i\}$ be the realized return; the better increases their bankroll by $\overline{a}_ix_i$ fraction if the wager wins, whereas otherwise they lose $x_i$ fraction.
The better knows the vector $\overline{\mathbf{a}}$ before placing their bets but does not know the probabilities of the different outcomes. 
Instead of relying on intuition to estimate the probabilities of outcomes for the different point spreads, we suppose that the better has access to a very accurate game emulator to run a finite number of game simulations to aid in deciding which wager or mix of wagers to bet on.

The goal of the better is to allocate their bankroll to maximize their overall expected return while also limiting their chance of losing money.
We pose this problem in the following way:
\begin{subequations}
    \label{eq:betting_prob}
    \begin{align}
        \max_\mathbf{x} \quad & \mathbb{E}\left[\mathbf{a}^\intercal\mathbf{x}\right] \\
        \text{s.t.} \quad & \text{Pr}\left(\mathbf{a}^\intercal\mathbf{x}\geq \beta \right) \geq 1-\alpha \\
        & \mathbf{1}^\intercal\mathbf{x} \leq 1 \\
        & \mathbf{x} \geq \mathbf{0}
    \end{align}
\end{subequations}
so that the probability of losing more than $\beta$ fraction of the bankroll is at most $\alpha$.

In our simulations, we simulate two independent games, each with two correlated wager options.
For each game $j\in\{1,2\}$, we use a uniform random variable $u_j\in[0,1]$ to represent the point outcome, and we use the probability of winning $\rho_{2(j-1)+i}$ to define the point spread of wager $i\in\{1,2\}$ so that:
\begin{align}
    a_{2(j-1)+i}(u_j) =
    \begin{cases}
    \overline{a}_i & \text{if } u_j\geq 1-\rho_{2(j-1)+i}, \\
    -1 & \text{otherwise.}
    \end{cases}
\end{align}
In other words, wager $i$ for game $j$ increases the bankroll if $u_j\geq 1-\rho_{2(j-1)+i}$.
The specific settings are $\boldsymbol{\rho}=[0.75,0.6,0.7,0.4]$, and $\overline{\mathbf{a}}=[0.5,0.95,0.6,2.1]$.

The metrics we use to evaluate the different methods are the average expected reward and the worst-case/maximum probability that $\mathbf{a}^\intercal\mathbf{x}<\beta$ among all trained solutions.
For a given sample size $N$, we draw $10^3$ different training sets, where each training set is processed and used to find a solution to Problem \eqref{eq:betting_prob} with a particular chance constraint method.
Each solution is evaluated with the same test set that has $10^6$ samples used to calculate its expected reward and its probability that $\mathbf{a}^\intercal\mathbf{x}<\beta$.
The chance constraint settings are $\alpha=0.2$ and $\beta=0.1$.
For methods that require the practitioner to choose a probability at which the chance constraint itself can be violated, we set it to 0.1 (indicated by parentheses in the figures).

\subsection{Results}

\begin{figure}
     \centering
     \begin{subfigure}[b]{0.6\textwidth}
         \centering
         \includegraphics[width=\textwidth]{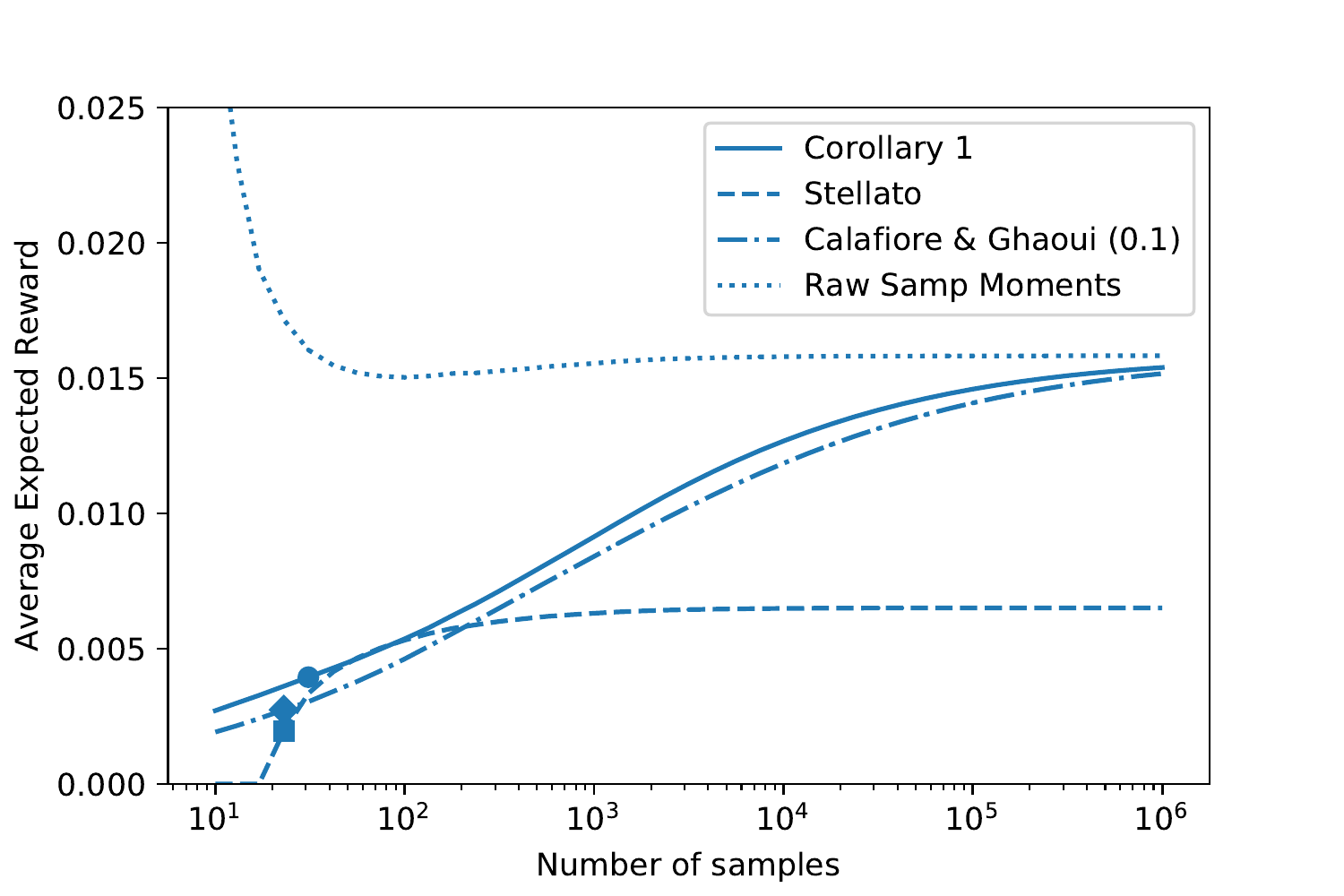}
         \caption{}
         \label{fig:Ereward_vs_N_cor1_stel_cala10_lem1}
     \end{subfigure}
     \hfill
     \begin{subfigure}[b]{0.6\textwidth}
         \centering
         \includegraphics[width=\textwidth]{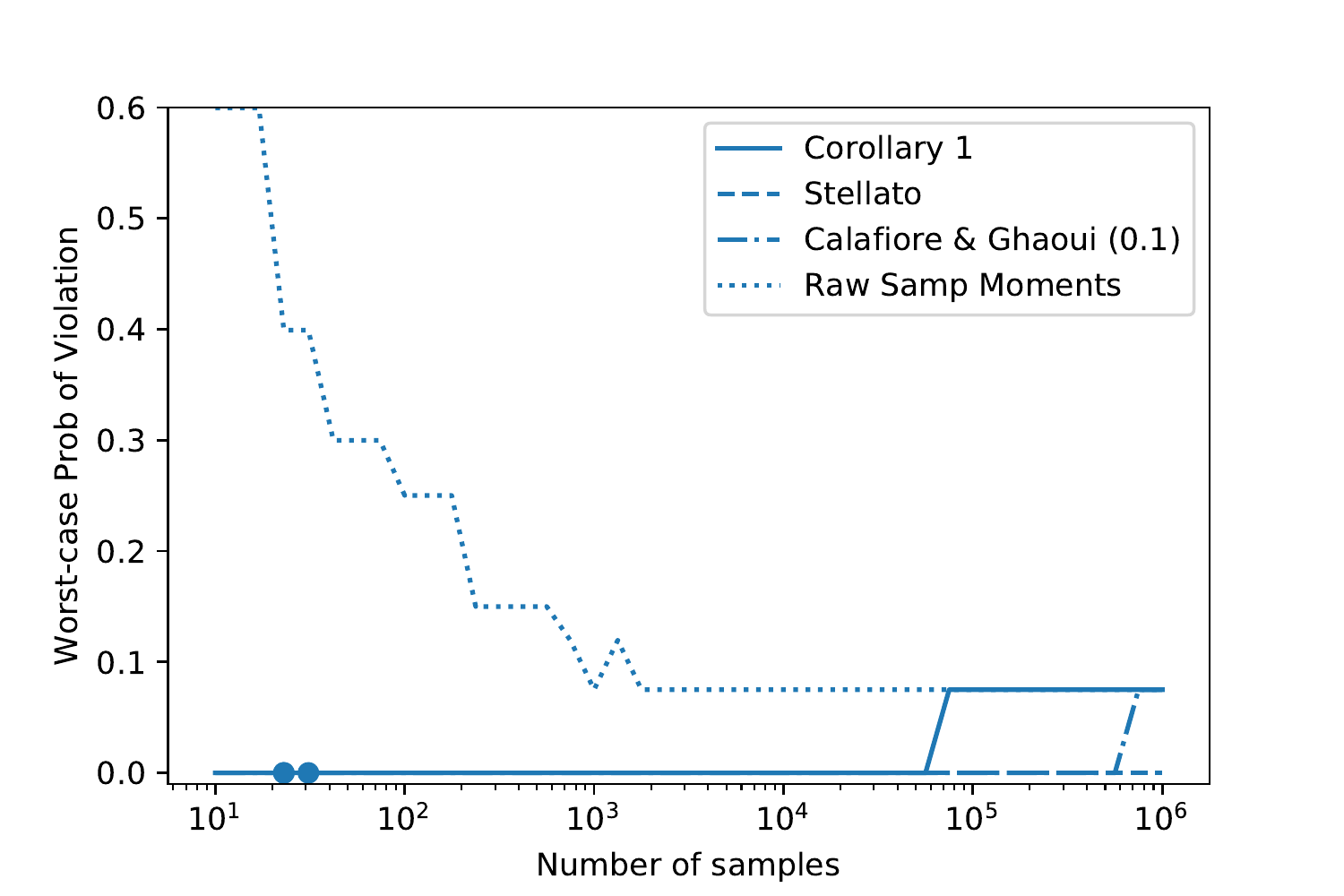}
         \caption{}
         \label{fig:maxprobvio_vs_N_cor1_stel_cala10}
     \end{subfigure}
    \caption{(a) Average expected reward and (b) worst-case probability that $\mathbf{a}^\intercal\mathbf{x}<\beta$ vs. number of samples.
    The dot/square/diamond indicates the minimum number of samples for the particular method's chance constraint guarantee. }
    \label{fig:momentbased}
\end{figure}

Figure \ref{fig:momentbased} shows both the average expected reward and the worst-case probability that $\mathbf{a}^\intercal\mathbf{x}<\beta$ versus the number of samples for different moment-based chance constraint methods.
It also includes marks that indicate the minimum number of samples needed for the associated theoretical guarantees to be applicable.
First, comparing our method set by Corollary \ref{thm:data_to_chance_const_p1} against simply applying the raw sample mean and covariance to the deterministic inequality in Lemma \ref{thm:calafiore3p1all}, we see that the average expected reward resulting from Corollary \ref{thm:data_to_chance_const_p1} asymptotically approaches that from Lemma \ref{thm:calafiore3p1all} as the number of samples increases, which corroborates the corollary.
Also, using the raw sample moments in Lemma \ref{thm:calafiore3p1all} allows for a possible solution to violate that chance constraint until the sample size is greater than 200, whereas Corollary \ref{thm:data_to_chance_const_p1} always satisfies the chance constraint.
Although using the raw sample moments might seem advantageous to Corollary \ref{thm:data_to_chance_const_p1} after a specific threshold number of samples, this threshold is not known a priori and is revealed only after running extensive simulations.

We also compare using Corollary \ref{thm:data_to_chance_const_p1} against a moment-based method close to ours developed by \cite{calafiore2006distributionally}, where a specific allowable probability that the chance constraint can be violated needs to be set. In our experiments, we set this probability to 0.1.
Although the method in \cite{calafiore2006distributionally} has asymptotic properties that are similar to Corollary \ref{thm:data_to_chance_const_p1} in these simulations, the chance constraint is only theoretically guaranteed up to one minus the set allowable probability; whereas Corollary \ref{thm:data_to_chance_const_p1} is guaranteed for it to be satisfied for any sample size greater than 26.

Not unlike our proposed approach, the Multivariate Sampled Chebyshev Approach by \cite{stellato2014data} does not require such a probability of violation to be set; in the results, however, the resulting deterministic constraint is so conservative that it asymptotically approaches one-third of the expected reward of Corollary \ref{thm:data_to_chance_const_p1}.

\begin{figure}
     \centering
     \begin{subfigure}[b]{0.6\textwidth}
         \centering
         \includegraphics[width=\textwidth]{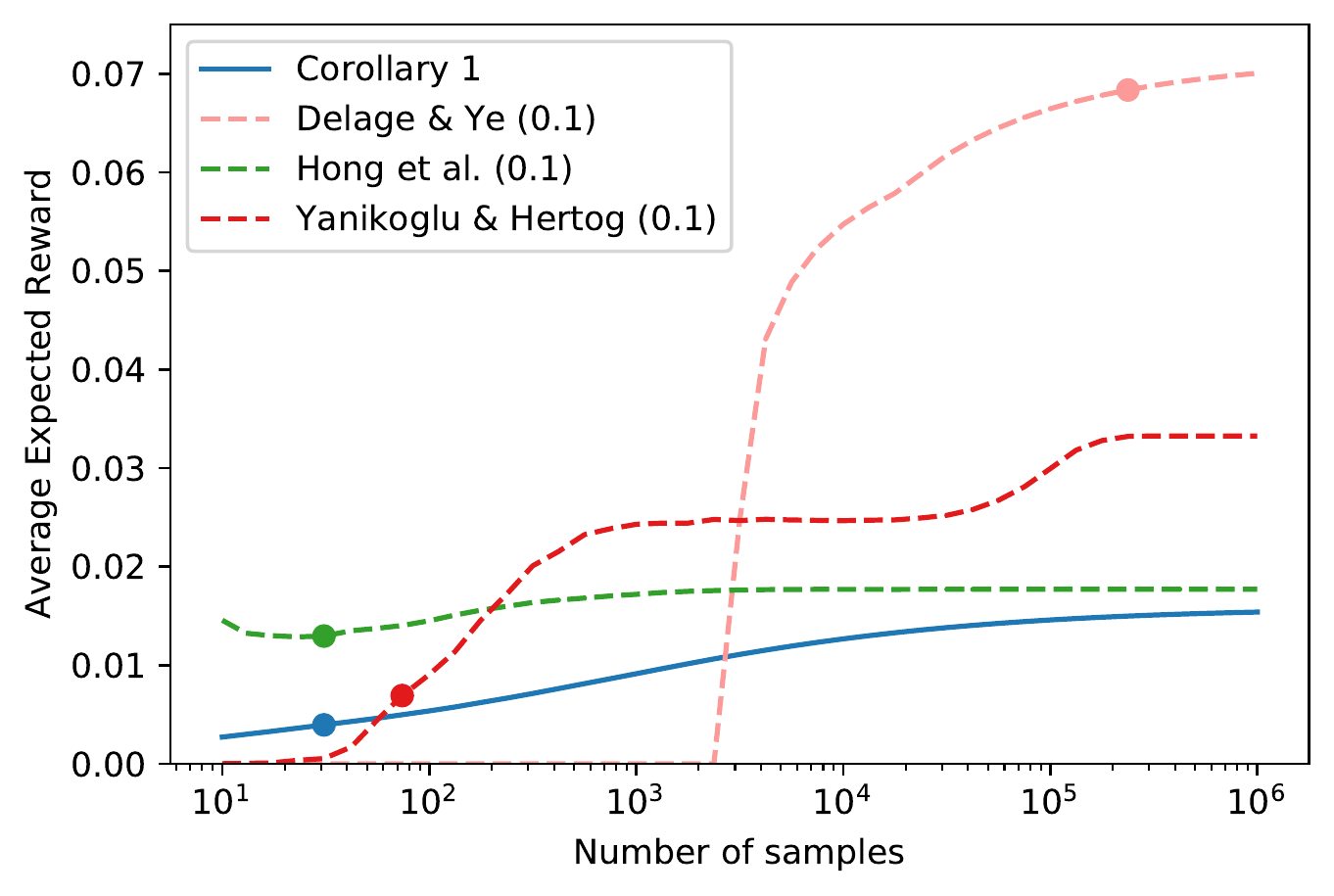}
         \caption{}
         \label{fig:Ereward_vs_N_cor1_dela_hong_yani}
     \end{subfigure}
     \hfill
     \begin{subfigure}[b]{0.6\textwidth}
         \centering
         \includegraphics[width=\textwidth]{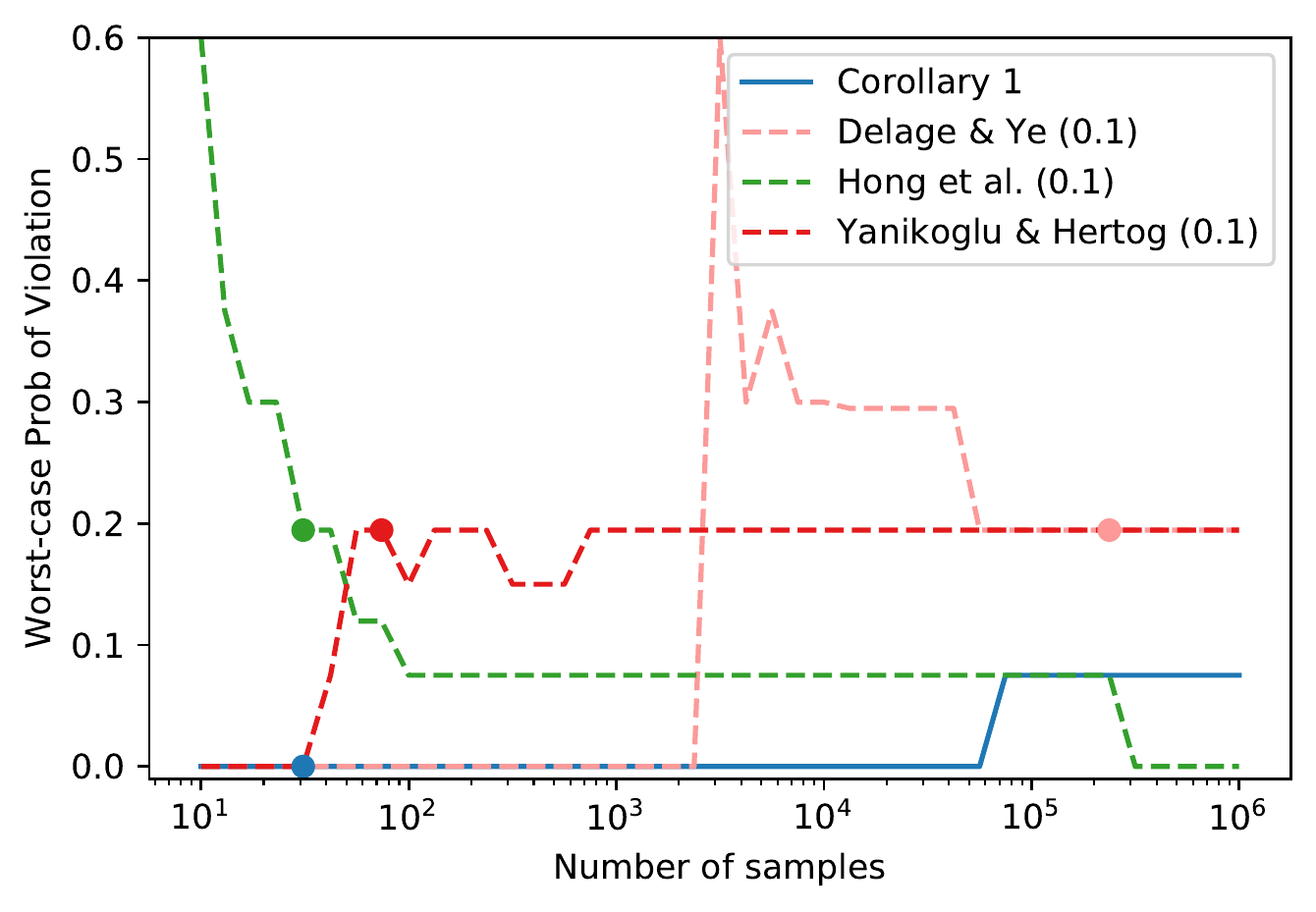}
         \caption{}
         \label{fig:maxprobvio_vs_N_cor1_dela_hong_yani}
     \end{subfigure}
    \caption{(a) Average expected reward and (b) worst-case probability that $\mathbf{a}^\intercal\mathbf{x}<\beta$ vs. number of samples.
    The dots indicate the minimum number of samples for the particular method's chance constraint guarantee. }
    \label{fig:cor1_dela_hong_yani}
\end{figure}

Additionally, in Figure \ref{fig:cor1_dela_hong_yani}, we compare Corollary \ref{thm:data_to_chance_const_p1} with more computationally expensive methods.
First, the moment uncertainty set method by \cite{delage2010distributionally} shows the most potential to having the tightest approximation because it has a much higher expected reward and it has a worst-case probability of violation that is only slightly less than $\alpha=0.2$ when the number of samples is greater than its minimum needed;
however, compared to Corollary \ref{thm:data_to_chance_const_p1}, it requires approximately $2000\times$ more samples for the method to validly apply without violating the chance constraint.
And because it adds $\theta(m^2)$ auxiliary variables to the optimization problem, it is much more computationally expensive to solve.

The union of discretized cells method by \cite{yanikouglu2013safe} also performs well and demonstrates a tight approximation of the chance constraint;
however, it requires the optimization problem to be solved $O(\sqrt{m)}$ times, where $m$ is the number of random variables, and it has the downside to requiring a problem-dependent way of defining the discretized cells.
In this problem, we simply defined a cell to be a specific realization of $a$ for which there were only 9 different options.
Had we increased the number of games or wagers, this would have exponentially increased the number of cells and thus increased the minimum number of samples needed because each cell needs at least five occurrences for the method's guarantees to be valid.
Of course, low-frequency cells could be merged to get the occurrences up to 5, but that adds extra computational overhead and additional tuning of the method.

The uncertainty set method by \cite{hong2017learning} uses the data to learn the size and shape of an uncertainty set.
One suggestion by them is to use the sample covariance matrix and mean to form an ellipsoid;
however, we found that using the identity matrix in place of the covariance matrix gave a less conservative approximation in this simulation.
Compared to Corollary \ref{thm:data_to_chance_const_p1}, it has only a slightly higher expected reward at high sample sizes but is significantly greater at lower sample sizes.
The reason behind this is that the uncertainty set method behaves like using raw sample moments in Lemma \ref{thm:calafiore3p1all} (See Figure \ref{fig:momentbased}) by starting out by violating the chance constraint in the worst case at low sample sizes but quickly satisfying them as the sample size increases.
Note that this is the opposite behavior of Corollary \ref{thm:data_to_chance_const_p1} and \cite{yanikouglu2013safe}, which starts off more conservatively and gradually allows more violation.

\begin{figure}
     \centering
     \begin{subfigure}[b]{0.6\textwidth}
         \centering
         \includegraphics[width=\textwidth]{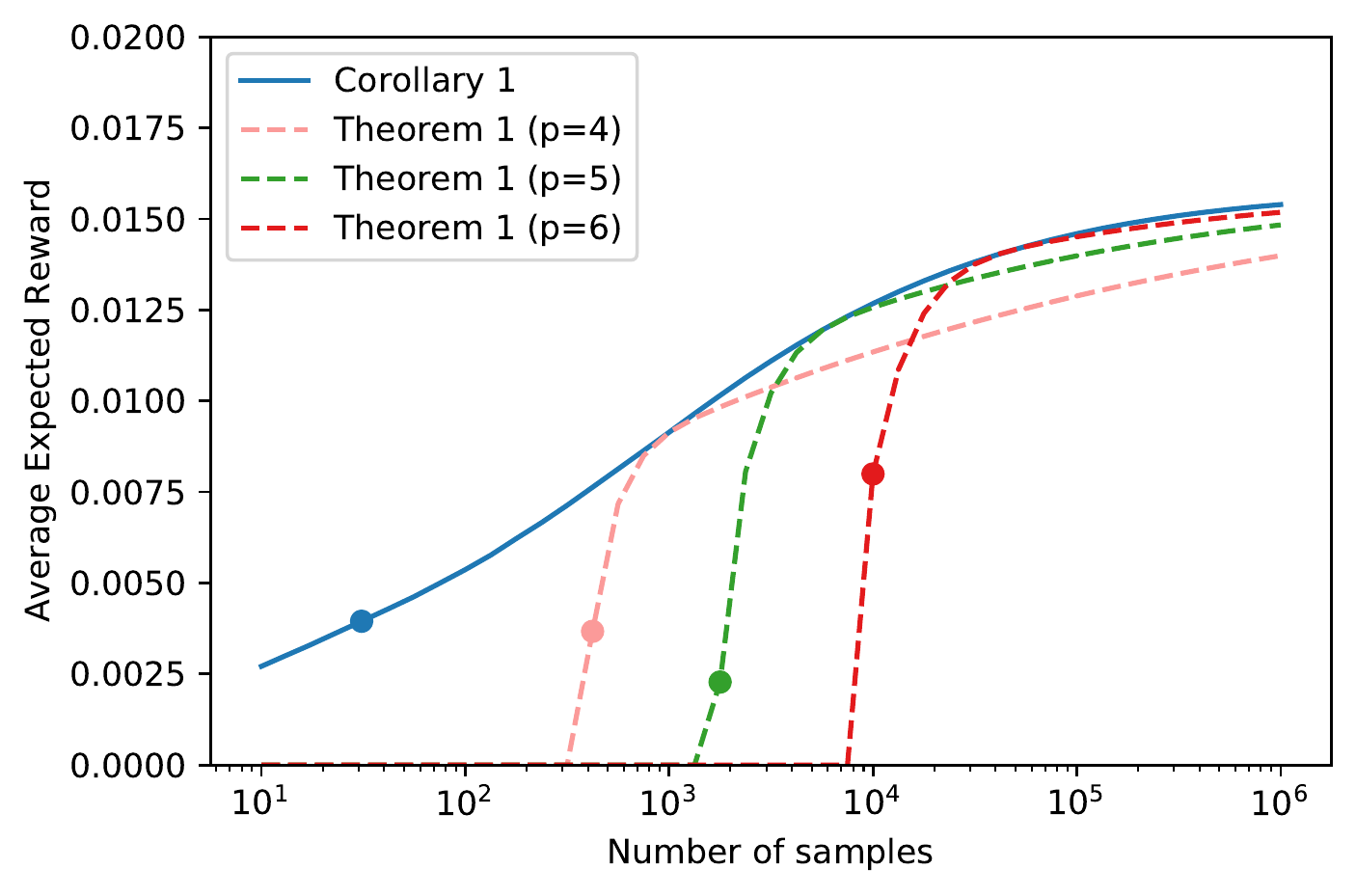}
         \caption{}
         \label{fig:Ereward_vs_N_cor1_thm1}
     \end{subfigure}
     \hfill
     \begin{subfigure}[b]{0.6\textwidth}
         \centering
         \includegraphics[width=\textwidth]{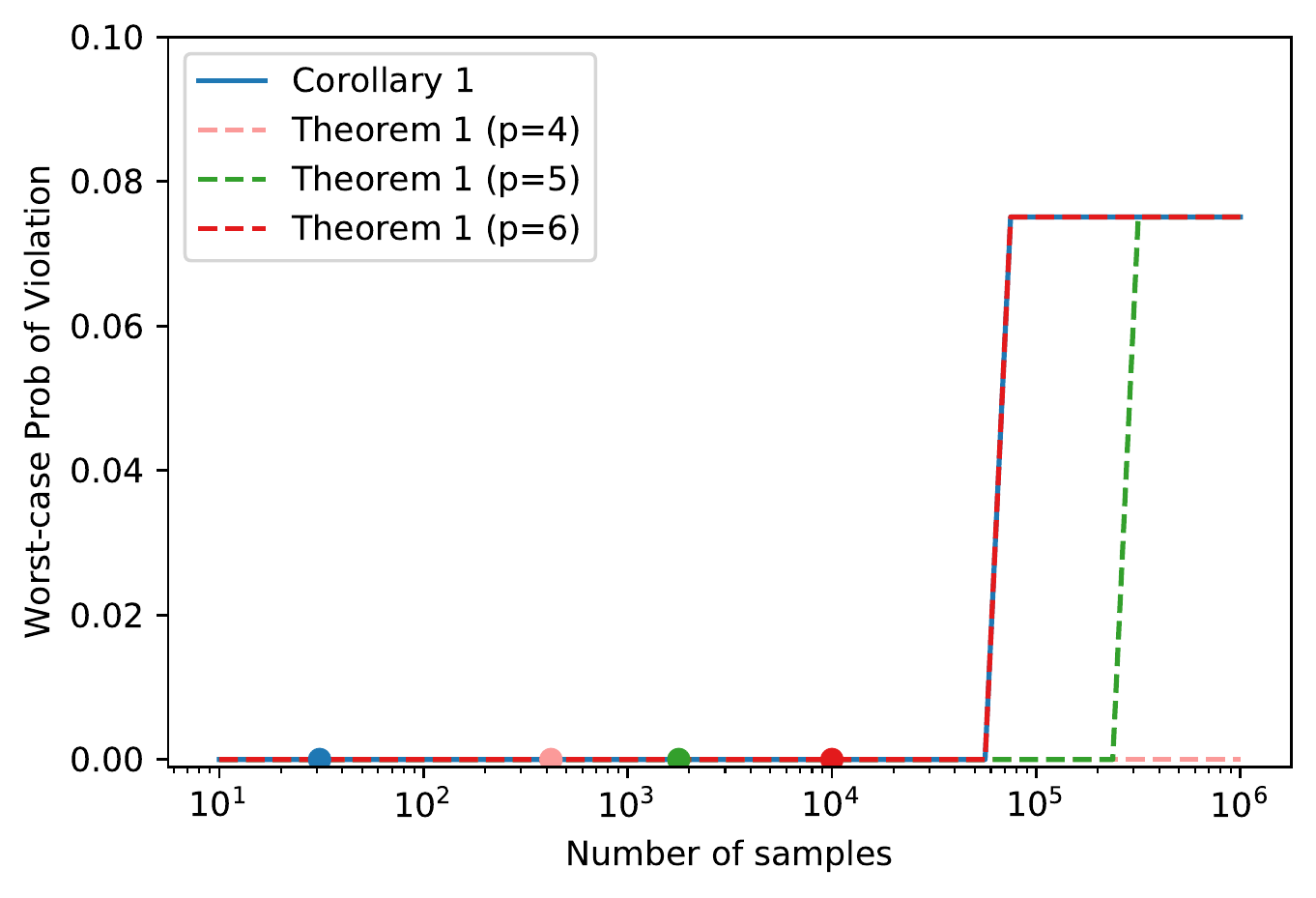}
         \caption{}
         \label{fig:maxprobvio_vs_N_cor1_thm1}
     \end{subfigure}
    \caption{(a) Average expected reward and (b) worst-case probability that $\mathbf{a}^\intercal\mathbf{x}<\beta$ vs. number of samples.}
    \label{fig:cor1_thm1}
\end{figure}

Figure \ref{fig:cor1_thm1} demonstrates the practical advantage of using Corollary \ref{thm:data_to_chance_const_p1} over setting a specific $p$ in Theorem \ref{thm:data_to_chance_const} besides not having to set $p$.
We see that as $p$ increases in Theorem \ref{thm:data_to_chance_const}, the average expected reward increases faster after a certain threshold number of samples;
however, the threshold at which it releases some of its conservativeness to give a nonzero average expected reward also increases.
Corollary \ref{thm:data_to_chance_const_p1} makes an upper envelope over the realizations of Theorem \ref{thm:data_to_chance_const} to get both the lowest sample size threshold and the fastest rate of increase on the average expected reward.
This same behavior is also seen with regards to the worst-case probability of violation.

\begin{figure}
    \centering
    \includegraphics[width=0.6\textwidth]{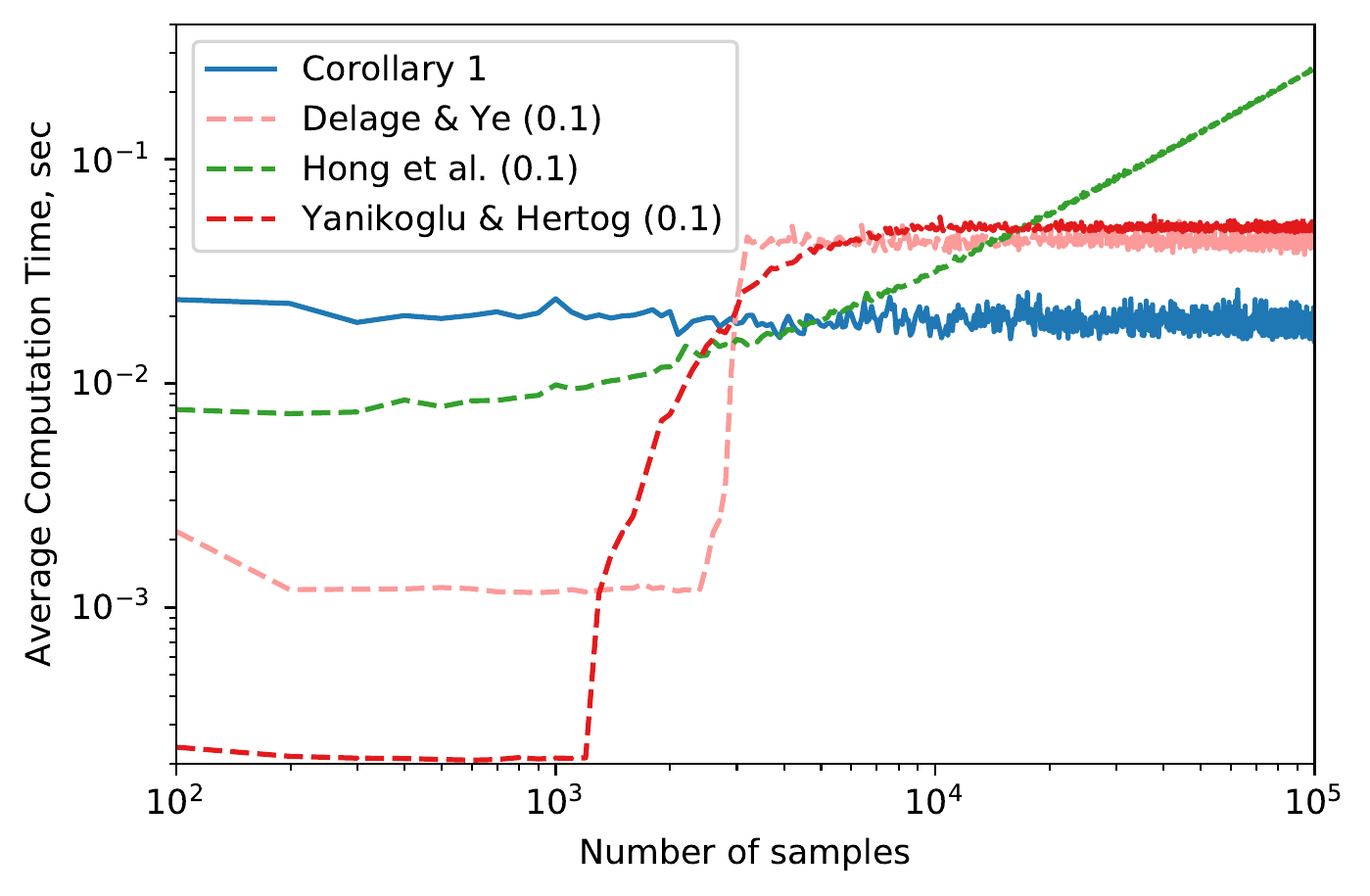}
    \caption{Average computation time versus number of acquired samples for different chance constraint methods.}
    \label{fig:seqRT}
\end{figure}

To demonstrate the computational time associated with sequential decision problems, we adjust our simulation's setup to emulate a sequential decision problem.
We assume that for every time step 100 new samples of $\mathbf{a}$ are recorded to aid in solving Problem \eqref{eq:betting_prob} with a particular chance constraint method.
This is run for 1,000 time steps, and thus $10^5$ samples are acquired.
At each time step, we measure the amount of time it takes to process the new data and solve for a new solution.
Figure \ref{fig:seqRT} gives the computation time averaged over 100 separate runs for the more computationally expensive methods.
First, we see that the computation time for Corollary \ref{thm:data_to_chance_const_p1} remains constant because it simply updates the sample mean and covariance with the 100 new samples and solves almost the same problems as the time step before.
On the other hand, the computational time of \cite{hong2017learning} increases with an increasing number of samples because all existing samples need to be reordered every time the sample mean is updated.
The computational time of \cite{delage2010distributionally} and \cite{yanikouglu2013safe} are very low at a low number of samples because there are not enough samples to trigger an optimization problem to be solved; instead, the methods simply implement the default safe option $\mathbf{x}=\mathbf{0}$, i.e., they bet none of the bankroll.
Once they have enough samples, however they solve an optimization problem---or mutliple optimization problems, in the case of \cite{yanikouglu2013safe}.
This makes the computational time almost double that of Corollary \ref{thm:data_to_chance_const_p1}.
In this small simulation example, where there are only four random variables, it might be a worthwhile trade-off to sacrifice computation time for more expected reward.
If the number of random variables $m$ is larger and $O(n)$, then the computation time (see Table \ref{tab:cc_complexity}) for solving the optimization problem becomes $O(n^4)$ for \cite{yanikouglu2013safe} and $O(n^7)$ for \cite{delage2010distributionally}, whereas the computation time is $O(n^{3.5})$ for the moment-based Corollary \ref{thm:data_to_chance_const_p1}.

\section{Conclusion}
\label{sec:conclusion}
In this paper, we designed a data-driven method for distributionally robust chance constraints that guarantees the satisfaction of the chance constraints when using the sample mean and covariance.
Along with the guarantees, we proved its asymptotic convergence as the number of samples approaches infinity to the distributionally robust chance constraint that has access to the true mean and covariance.
Although most other data-driven methods guarantee satisfaction of the chance constraint only up to a user-specified probability, our proposed method guarantees satisfaction of the chance constraint after a small number of samples are obtained.
We evaluated our data-driven method on a numerical example and showed that it (a) guarantees the satisfaction of the chance constraints, (b) approaches the optimal solution of the case when the true mean and covariance are known, and (c) has a computational advantage compared with other state-of-the-art methods.
Future work can take our data-driven method in several directions.
For example, our method could be adjusted to the case when the mean and covariance are time-varying, which would broaden the usefulness of the application.

\appendix
\section{Proof for Radius of the Estimated Support of an Ellipsoid}
\label{prf:r_ellip}

\begin{proposition}
    \label{thm:r_ellip}
    If the estimated support $\hat{\mathcal{S}}$ of the random variable $\mathbf{a}$ is defined by the ellipsoid $(\mathbf{a}-\mathbf{c})^\intercal\mathbf{V}(\mathbf{a}-\mathbf{c})\leq 1$, where $\mathbf{c}$ is its center, and $\mathbf{V}$ is a symmetric positive definite matrix, then the function of the radius of $\hat{\mathcal{S}}$ with respect to the weighting $\mathbf{x}$, i.e., $r(\mathbf{x}):=\frac{1}{2}\sup_{\{\mathbf{a}_1,\mathbf{a}_2\}\in\hat{\mathcal{S}}}|\mathbf{a}_1^\intercal\mathbf{x}-\mathbf{a}_2^\intercal\mathbf{x}|$, is equal to $\sqrt{\mathbf{x}^\intercal\mathbf{V}^{-1}\mathbf{x}}$.
\end{proposition}

    First, the matrix $\mathbf{V}$ is invertible because it is symmetric positive definite matrix.
    Note that for any $\mathbf{a}^{'}$ and $\mathbf{a}^{''}$ in $\hat{\mathcal{S}}$, one of $(\mathbf{a}^{'} - \mathbf{a}^{''})^\intercal \mathbf{x}$ and $(\mathbf{a}^{''} - \mathbf{a}^{'})^\intercal \mathbf{x}$ is nonnegative. Therefore, without loss of generality, $r(\mathbf{x})$ can be equivalently written in the following linear optimization problem because of the problem's symmetry and the closedness of $\hat{\mathcal{S}}$:
    \begin{subequations}
        \begin{align}
            r(\mathbf{x}) := \max_{\mathbf{a}_1,\mathbf{a_2}}~ & \frac{1}{2}\mathbf{x}^\intercal (\mathbf{a}_1-\mathbf{a}_2) \nonumber \\
            \text{s.t. }~
            & (\mathbf{a}_1-\mathbf{c})^\intercal\mathbf{V}(\mathbf{a}_1-\mathbf{c}) \leq 1 \nonumber \\
            & (\mathbf{a}_2-\mathbf{c})^\intercal\mathbf{V}(\mathbf{a}_2-\mathbf{c}) \leq 1. \nonumber
        \end{align}
    \end{subequations}
    We can use the Karush–Kuhn–Tucker (KKT) conditions and take advantage that the problem is separable between $\mathbf{a}_1$ and $\mathbf{a}_2$ to find the optimal value.
    
    The first problem is:
    \begin{subequations}
        \begin{align}
            \max_{\mathbf{a}_1}~ & \frac{1}{2}\mathbf{x}^\intercal \mathbf{a}_1 \nonumber \\
            \text{s.t. }~
            & (\mathbf{a}_1-\mathbf{c})^\intercal\mathbf{V}(\mathbf{a}_1-\mathbf{c}) \leq 1 \nonumber
        \end{align}
    \end{subequations}
    which gives the KKT conditions:
    \begin{subequations}
        \begin{align}
            \frac{1}{2}\mathbf{x}^\intercal - 2\lambda_1(\mathbf{a}_1-\mathbf{c})^\intercal\mathbf{V} & = 0 \nonumber \\
            \lambda_1 & \geq 0 \nonumber \\
            \lambda_1((\mathbf{a}_1-\mathbf{c})^\intercal\mathbf{V}(\mathbf{a}_1-\mathbf{c}) - 1) & = 0 \nonumber \\
            (\mathbf{a}_1-\mathbf{c})^\intercal\mathbf{V}(\mathbf{a}_1-\mathbf{c}) - 1 & \leq 0 \nonumber
        \end{align}
    \end{subequations}
    and the unique solution for $\mathbf{x}\neq \mathbf{0}$:
    \begin{subequations}
        \begin{align}
            \lambda_1 & = \frac{1}{4}\sqrt{\mathbf{x}^\intercal\mathbf{V}^{-1}\mathbf{x}} \nonumber \\
            \mathbf{a}_1 & = \mathbf{c} + \frac{1}{\sqrt{\mathbf{x}^\intercal\mathbf{V}^{-1}\mathbf{x}}}\mathbf{V}^{-1}\mathbf{x}, \nonumber
        \end{align}
    \end{subequations}
    and $\lambda_1=0$, with $\mathbf{a}_1=\mathbf{c}$ as a (nonunique) solution when $\mathbf{x}=\mathbf{0}$.
    
    The second problem is:
    \begin{subequations}
        \begin{align}
            \min_{\mathbf{a}_2}~ & \frac{1}{2}\mathbf{x}^\intercal \mathbf{a}_2 \nonumber \\
            \text{s.t. }~
            & (\mathbf{a}_2-\mathbf{c})^\intercal\mathbf{V}(\mathbf{a}_2-\mathbf{c}) \leq 1 \nonumber
        \end{align}
    \end{subequations}
    which gives the KKT conditions:
    \begin{subequations}
        \begin{align}
            \frac{1}{2}\mathbf{x}^\intercal + 2\lambda_2(\mathbf{a}_2-\mathbf{c})^\intercal\mathbf{V} & = 0 \nonumber \\
            \lambda_2 & \geq 0 \nonumber \\
            \lambda_2((\mathbf{a}_2-\mathbf{c})^\intercal\mathbf{V}(\mathbf{a}_2-\mathbf{c}) - 1) & = 0 \nonumber \\
            (\mathbf{a}_2-\mathbf{c})^\intercal\mathbf{V}(\mathbf{a}_2-\mathbf{c}) - 1 & \leq 0 \nonumber
        \end{align}
    \end{subequations}
    and the unique solution for $\mathbf{x}\neq \mathbf{0}$:
    \begin{subequations}
        \begin{align}
            \lambda_2 & = \frac{1}{4}\sqrt{\mathbf{x}^\intercal\mathbf{V}^{-1}\mathbf{x}} \nonumber \\
            \mathbf{a}_2 & = \mathbf{c} - \frac{1}{\sqrt{\mathbf{x}^\intercal\mathbf{V}^{-1}\mathbf{x}}}\mathbf{V}^{-1}\mathbf{x}, \nonumber
        \end{align}
    \end{subequations}
    and $\lambda_2=0$  with $\mathbf{a}_2=\mathbf{c}$ as a (nonunique) solution when $\mathbf{x}=\mathbf{0}$.
    
    Putting both solutions together, we get the radius to be a convex function of $\mathbf{x}$:
    \begin{align}
        r(\mathbf{x}) & = \frac{1}{2}\mathbf{x}^\intercal (\mathbf{a}_1-\mathbf{a}_2) \nonumber \\
        & = \sqrt{\mathbf{x}^\intercal\mathbf{V}^{-1}\mathbf{x}}. \nonumber
    \end{align}
    \hfill $\qed$

\section{Proof of Theorem \ref{thm:data_to_chance_const}}
\label{prf:data_to_chance_const}

    First, define the random variable $v:=\mathbf{a}^\intercal\mathbf{x}$ coming from $\mathbf{a}$, which then has mean $\boldsymbol{\mu}^{*\intercal}\mathbf{x}$ denoted by $\beta^*$ and variance $\mathbf{x}^\intercal\boldsymbol{\Sigma}^*\mathbf{x}$ denoted by $\gamma^*$.
    From $N$ independent samples of $\mathbf{a}$, we can calculate the associated sample mean and the variance of $v$, which are $\hat{\beta}_N=\hat{\boldsymbol{\mu}}_N^\intercal\mathbf{x}$ and $\hat{\gamma}_N=\mathbf{x}^\intercal\hat{\boldsymbol{\Sigma}}_N\mathbf{x}$, respectively.
    Then Equation \eqref{eq:prob_cc_distrob_true} can be rewritten in its univariate form:
    \begin{align}
        \inf_{P\in\mathcal{P}(\beta^*,\gamma^*)}\text{Pr}\left(v\leq 0|v\sim P\right) \geq 1- \alpha.
        \label{eq:proof2_DRCC_prob_perf_uni}
    \end{align}
    Let $\mathcal{Q}_{\delta,N}$ be a confidence region of the first two moments, such that $(\beta^*,\gamma^*)\in\mathcal{Q}_{\delta,N}$, with a probability of at least $1-\delta$ for $\delta\in(0,1)$.
    Then by setting $c:=0$, and $\epsilon:=\frac{\alpha-\delta}{1-\delta}$ for $\delta\in(0,\alpha)$ in Lemma \ref{thm:confregion}, we have that if
    \begin{align}
        \inf_{P\in\mathcal{P}(\beta,\gamma)}\text{Pr}\left(v \leq 0 | v\sim P \right) \geq 1-\frac{\alpha-\delta}{1-\delta}, \quad \forall (\beta,\gamma)\in \mathcal{Q}_{\delta,N} \label{eq:proof2_DRCC_prob_conf}
    \end{align}
    is satisfied, then so is Equation \eqref{eq:proof2_DRCC_prob_perf_uni}.
    By setting $\mathbf{x}$ to the scalar 1, $\boldsymbol{\mu}^*$ to the scalar $\beta$, and  $\boldsymbol{\Sigma}^*$ to the scalar $\gamma$ in Lemma \ref{thm:calafiore3p1all}, we have that
    \begin{align}
        \beta + \sqrt{\frac{1-\alpha}{\alpha-\delta}}\sqrt{\gamma} \leq 0, \quad \forall (\beta,\gamma)\in \mathcal{Q}_{\delta,N}
    \end{align}
    is equivalent to Equation \eqref{eq:proof2_DRCC_prob_conf}.
    Letting $\mathcal{Q}_{\delta,N}$ be a closed set, the above inequality is only satisfied for every case in $\mathcal{Q}_{\delta,N}$ if:
    \begin{align}
        \max_{(\beta,\gamma)\in \mathcal{Q}_{\delta,N}}\left\{\beta + \sqrt{\frac{1-\alpha}{\alpha-\delta}}\sqrt{\gamma}\right\} \leq 0. \label{eq:proof2_maxQ}
    \end{align}
    Using Lemma \ref{thm:calafiore4p1} with $v$ as the random variable and the assumption that $\mathcal{S}^*\subseteq\hat{\mathcal{S}}$, which makes $r(\mathbf{x})\geq r^*$, we can define a specific form for the confidence region around $(\hat{\beta}_N,\hat{\gamma}_N)$:
    \begin{align}
        \mathcal{Q}_{\delta,N} := & \Big\{(\beta,\gamma):|\beta-\hat{\beta}_N| \leq \frac{r(\mathbf{x})}{\sqrt{N}}\left(2+\sqrt{2\ln(2/\delta)}\right), \nonumber \\ & \quad |\gamma-\hat{\gamma}_N| \leq \frac{2(r(\mathbf{x}))^2}{\sqrt{N}}\left(2+\sqrt{2\ln(4/\delta)}\right)\Big\}. \nonumber
    \end{align}
    Then we can maximize the first and last terms on the LHS of \eqref{eq:proof2_maxQ} because the structure of $\mathcal{Q}_{\delta,N}$ allows both terms to be maximized independently to get:
    \begin{align}
        \hat{\beta}_N + \frac{r(\mathbf{x})}{\sqrt{N}}\left(2+\sqrt{2\ln(2/\delta)}\right) \quad\quad\quad\quad\quad\quad\quad\quad\quad\quad\quad\quad\quad\quad\quad\quad \nonumber \\
        + \sqrt{\frac{1-\alpha}{\alpha-\delta}}\sqrt{\hat{\gamma}_N+\frac{2(r(\mathbf{x}))^2}{\sqrt{N}}\left(2+\sqrt{2\ln(4/\delta)}\right)}\leq 0. \label{eq:proof2_determin_samp}
    \end{align}
    Next, we increase the second factor of the second term from $(2+\sqrt{2\ln(2/\delta)})$ to $(2+\sqrt{2\ln(4/\delta)})$, which makes the constraint more conservative but simplifies the analysis in the next steps, and we set $\delta:=4\exp\left(-(N^\frac{1}{p}-2)^2/2\right)$.
    This makes Equation \eqref{eq:proof2_determin_samp} become:
    \begin{align}
        \hat{\boldsymbol{\mu}}_N^\intercal\mathbf{x} + \phi_N r(\mathbf{x}) + \kappa_N\sqrt{\frac{1-\alpha}{\alpha}}\sqrt{\mathbf{x}^\intercal\hat{\boldsymbol{\Sigma}}_N\mathbf{x}+2\phi_N (r(\mathbf{x}))^2} \leq 0, \label{eq:proof2_DRCC_det}
    \end{align}
    and is thus a conservative approximation of \eqref{eq:prob_cc_distrob_true}.
    Note that by using this specific setting of $\delta$, the condition on the number of samples in Lemma \ref{thm:calafiore4p1} is satisfied because $\left(2+\sqrt{2\ln(4/\delta)}\right)^2=N^\frac{2}{p}\leq N$ for $p>2$ and $N\geq 1$.
    The condition on $N$ in the theorem's statement makes sure that $\delta<\alpha$.
    
    As $N\rightarrow\infty$, then $\hat{\boldsymbol{\mu}}\rightarrow\boldsymbol{\mu}^*$ and $\hat{\boldsymbol{\Sigma}}\rightarrow\boldsymbol{\Sigma}^*$ because of the law of large numbers.
    This results in that as $N\rightarrow\infty$, then $\kappa_N\rightarrow 1$, and Equation \eqref{eq:proof2_DRCC_det} asymptotically approaches:
    \begin{align}
        \boldsymbol{\mu}^{*\intercal}\mathbf{x} + \sqrt{\frac{1-\alpha}{\alpha}}\sqrt{\mathbf{x}^\intercal\boldsymbol{\Sigma}^*\mathbf{x}} \leq 0 \label{eq:proof2_perf_det}
    \end{align}
    which is equivalent to \eqref{eq:prob_cc_distrob_true} according to Lemma \ref{thm:calafiore3p1all}.
    
    Finally, the transformation of Equation \eqref{eq:proof2_DRCC_det} to \eqref{eq:prob_cc_distrob_data} comes from replacing  $\mathbf{x}^\intercal\hat{\boldsymbol{\Sigma}}_N\mathbf{x}$ with $y_1^2$ and $2\phi_N r(\mathbf{x})^2$ with $y_2^2$ in Equation \eqref{eq:proof2_DRCC_det}, and then using $\mathbf{y}$ as upper bounds.
    This makes constraints \eqref{eq:prob_cc_distrob_data} convex respect to $\mathbf{x}$, and note that $\mathbf{r}(\mathbf{x})$ is convex because of its definition (See \cite{boyd2004convex}, Example 3.7).
    
    \hfill $\qed$

\section{Proof of Lemma \ref{thm:confregion}}
\label{prf:confregion}

    The confidence region and its probabilistic relationship is formally stated as:
    \begin{align}
        \text{Pr}\left(\mathcal{Q}_{\delta}\ni(\mu^*,\sigma^{2*}) \right) \geq 1-\delta. \label{eq:proof2_conf_reg_prob}
    \end{align}
    Equation \eqref{eq:proof2_const_DRO_samp} means that for any distribution $P$ in $\mathcal{P}(\mu,\sigma^2,\underline{v}^*,\overline{v}^*)$, where $(\mu,\sigma^2)\in\mathcal{Q}_\delta$, we have that:
    \begin{align}
        \text{Pr}\left(v + c\leq 0 | v\sim P\in\mathcal{P}(\mu,\sigma^2,\underline{v}^*,\overline{v}^*) \right) \geq 1-\epsilon, \quad \forall (\mu,\sigma^2)\in \mathcal{Q}_\delta \nonumber
    \end{align}
    because of the infimum.
    Thus, for any case when $\mathcal{Q}_{\delta}$ contains $(\mu,\sigma^2)$, we consequently have that:
    \begin{align}
        \text{Pr}\left(v+c \leq 0 | (v\sim P\in\mathcal{P}(\mu,\sigma^2,\underline{v}^*,\overline{v}^*)) \cap (\mathcal{Q}_{\delta}\ni(\mu,\sigma^2))  \right) \geq 1-\epsilon. \label{eq:proof2_const_DRO_conf}
    \end{align}
    Next, we bound the probability of violation for any distribution $P$ in $\mathcal{P}(\mu^*,\sigma^{2*},\underline{v}^*,\overline{v}^*)$, and then we take the worst-case via the infimum at the end.
    With the law of total probability, we can condition on whether the confidence region $\mathcal{Q}_{\delta}$ contains $(\mu^*,\sigma^{2*})$ or not, and then we apply the above bounds:
    \begin{align}
        \text{Pr}\left(v+c \leq 0 | v\sim P\in\mathcal{P}(\mu^*,\sigma^{2*},\underline{v}^*,\overline{v}^*)) \right) \quad\quad\quad\quad\quad\quad\quad\quad\quad\quad\quad\quad\quad\quad\quad \nonumber
    \end{align}
    \begin{align}
        & = \text{Pr}\left(v+c \leq 0 | (v\sim P\in\mathcal{P}(\mu^*,\sigma^{2*},\underline{v}^*,\overline{v}^*)) \cap (\mathcal{Q}_{\delta}\ni(\mu^*,\sigma^{2*})) \right) \nonumber \\
        & \quad\quad\quad\quad \text{Pr}\left(\mathcal{Q}_{\delta}\ni(\mu^*,\sigma^{2*}) | v\sim P\in\mathcal{P}(\mu^*,\sigma^{2*},\underline{v}^*,\overline{v}^*) \right) \nonumber \\
        & \quad\quad + \text{Pr}\left(v+c \leq 0 | (v\sim P\in\mathcal{P}(\mu^*,\sigma^{2*},\underline{v}^*,\overline{v}^*)) \cap (\mathcal{Q}_{\delta}\not\ni(\mu^*,\sigma^{2*})) \right) \nonumber \\
        & \quad\quad\quad\quad \text{Pr}\left(\mathcal{Q}_{\delta}\not\ni(\mu^*,\sigma^{2*}) | v \sim P\in\mathcal{P}(\mu^*,\sigma^{2*},\underline{v}^*,\overline{v}^*) \right) \nonumber \\
        & \geq \text{Pr}\left(v+c \leq 0 | (v\sim P\in\mathcal{P}(\mu^*,\sigma^{2*},\underline{v}^*,\overline{v}^*)) \cap (\mathcal{Q}_{\delta}\ni(\mu^*,\sigma^{2*})) \right) \nonumber \\
        & \quad\quad\quad\quad \text{Pr}\left(\mathcal{Q}_{\delta}\ni(\mu^*,\sigma^{2*}) | v\sim P\in\mathcal{P}(\mu^*,\sigma^{2*},\underline{v}^*,\overline{v}^*) \right) \nonumber \\
        & = \text{Pr}\left(v+c \leq 0 | (v\sim P\in\mathcal{P}(\mu^*,\sigma^{2*},\underline{v}^*,\overline{v}^*)) \cap (\mathcal{Q}_{\delta}\ni(\mu^*,\sigma^{2*})) \right) \nonumber \\
        & \quad\quad\quad\quad \text{Pr}\left(\mathcal{Q}_{\delta}\ni(\mu^*,\sigma^{2*})\right) \nonumber \\
        & \geq (1-\epsilon)\text{Pr}\left(\mathcal{Q}_{\delta}\ni(\mu^*,\sigma^{2*})\right) \nonumber \\
        & \geq (1-\epsilon)(1-\delta).
    \end{align}
    The second inequality applies Equation \eqref{eq:proof2_const_DRO_conf}, and the third inequality applies Equation \eqref{eq:proof2_conf_reg_prob}.
    Because the above inequality is true for any distribution $P$ in the family with mean $\mu^*$, variance $\sigma^{2*}$, and support $[\underline{v}^*,\overline{v}^*]$, then we can say that the infimum over all the distributions is Equation \eqref{eq:proof2_const_DRO_perf}.
    
    \hfill $\qed$

\section{Proof of Proposition \ref{thm:data_to_chance_const_ind}}
\label{prf:data_to_chance_const_ind}

    First, define the random variable $v:=\mathbf{a}^\intercal\mathbf{x}$ coming from $\mathbf{a}$, which then has mean $\boldsymbol{\mu}^{*\intercal}\mathbf{x}$ denoted by $\beta^*$, and supports $\underline{v}^*:=\sum_i^n\min\{\underline{s}_i^*x_i,\overline{s}_i^*x_i\}$ and $\overline{v}^*:=\sum_i^n\max\{\underline{s}_i^*x_i,\overline{s}_i^*x_i\}$.
    From $N$ independent samples of $\mathbf{a}$, we can calculate the associated sample mean of $v$, which is $\hat{\beta}_N=\hat{\boldsymbol{\mu}}_N^\intercal\mathbf{x}$.
    Then Equation \eqref{eq:prob_cc_distrob_true_ind} can be rewritten in its univariate form:
    \begin{align}
        \inf_{P\in\mathcal{P}(\beta^*,\underline{v}^*,\overline{v}^*)}\text{Pr}\left(v\leq 0|v\sim P\right) \geq 1- \alpha
        \label{eq:proof2_DRCC_prob_perf_ind_uni}
    \end{align}
    where the variance is left unbounded.
    Let $[\underline{q},\overline{q}]_{\delta,N}$ be a confidence interval of the mean such that $\beta^*\in[\underline{q},\overline{q}]_{\delta,N}$ with a probability of at least $1-\delta$ for $\delta\in(0,1)$.
    Then by setting $c:=0$, $\epsilon:=\frac{\alpha-\delta}{1-\delta}$ for $\delta\in(0,\alpha)$, and $\mathcal{Q}_\delta:=\{(\mu,\sigma^2)\in\mathbb{R}^2:\mu\in[\underline{q},\overline{q}]_{\delta,N}\}$ in Lemma \ref{thm:confregion}, we have that if
    \begin{align}
        \inf_{P\in\mathcal{P}(\beta,\underline{v}^*,\overline{v}^*)}\text{Pr}\left(v \leq 0 | v\sim P \right) \geq 1-\frac{\alpha-\delta}{1-\delta}, \quad \forall \beta\in [\underline{q},\overline{q}]_{\delta,N} \label{eq:proof2_DRCC_prob_ind_conf}
    \end{align}
    is satisfied, then so is Equation \eqref{eq:proof2_DRCC_prob_perf_ind_uni}.
    By setting $\mathbf{x}$ to the scalar 1, $\boldsymbol{\mu}^*$ to the scalar $\beta$, and  $\mathbf{S}^*$ to the scalar $\overline{v}^*-\underline{v}^*$ in Lemma \ref{thm:calafiore3p2}, we have that if
    \begin{align}
        \beta  + \sqrt{\frac{1}{2}\ln\left(\frac{1-\delta}{\alpha-\delta}\right)}(\overline{v}^*-\underline{v}^*) \leq 0, \quad \forall \beta\in [\underline{q},\overline{q}]_{\delta,N}
    \end{align}
    is satisfied, then so is Equation \eqref{eq:proof2_DRCC_prob_ind_conf}.
    Because $[\underline{q},\overline{q}]_{\delta,N}$ is a closed set, the above inequality is only satisfied for every case in $[\underline{q},\overline{q}]_{\delta,N}$ if
    \begin{align}
        \max_{\beta\in [\underline{q},\overline{q}]_{\delta,N}}\left\{\beta +\sqrt{\frac{1}{2}\ln\left(\frac{1-\delta}{\alpha-\delta}\right)}(\overline{v}^*-\underline{v}^*)\right\} \leq 0. \label{eq:proof2_maxq_ind}
    \end{align}
    Using Lemma \ref{thm:shawe} with $v$ as the random variable and the assumption that $\underline{\mathbf{s}}\leq\underline{\mathbf{s}}^*$ and $\overline{\mathbf{s}}\geq\overline{\mathbf{s}}^*$, which makes $r_\text{int}(\mathbf{x}) = \frac{1}{2}\|\mathbf{S}\mathbf{x}\|_1\geq r^*$, we can define a specific form for the confidence interval around $\hat{\beta}_N$:
    \begin{align}
        [\underline{q},\overline{q}]_{\delta,N} := & \Big\{\beta:|\beta-\hat{\beta}_N| \leq \frac{\frac{1}{2}\|\mathbf{S}\mathbf{x}\|_1}{\sqrt{N}}\left(2+\sqrt{2\ln(1/\delta)}\right)\Big\}. \nonumber
    \end{align}
    Then we can maximize the first term on the LHS of \eqref{eq:proof2_maxq_ind} because the structure of $[\underline{q},\overline{q}]_{\delta,N}$ is one-dimensional, it allows the term to easily be maximized to get:
    \begin{align}
        \hat{\beta}_N + \frac{\frac{1}{2}\|\mathbf{S}\mathbf{x}\|_1}{\sqrt{N}}\left(2+\sqrt{2\ln(1/\delta)}\right) + \sqrt{\frac{1}{2}\ln\left(\frac{1-\delta}{\alpha-\delta}\right)}\|\mathbf{S}\mathbf{x}\|_1  \leq 0 \label{eq:proof2_determin_ind_samp}
    \end{align}
    where the third term above conservatively replaced $\overline{v}^*-\underline{v}^*$ with $\|\mathbf{S}\mathbf{x}\|_1$ because:
    \begin{align}
        \overline{v}^*-\underline{v}^* & = \sum_i^n(\max\{\underline{s}_i^*x_i,\overline{s}_i^*x_i\} - \min\{\underline{s}_i^*x_i,\overline{s}_i^*x_i\}) \nonumber \\
        & = \sum_i^n\max\{(\underline{s}_i^*-\overline{s}_i^*)x_i,(\overline{s}_i^*-\underline{s}_i^*)x_i\}  \nonumber \\
        & = \sum_i^n|(\overline{s}_i^*-\underline{s}_i^*)x_i| \nonumber \\
        & \leq \sum_i^n|(\overline{s}_i-\underline{s}_i)x_i| = \|\mathbf{S}\mathbf{x}\|_1. \nonumber
    \end{align}
    Next, we set $\delta:=\exp\left(-(N^\frac{1}{p}-2)^2/2\right)$.
    This makes Equation \eqref{eq:proof2_determin_ind_samp} become \eqref{eq:prob_cc_distrob_data_ind}
    and is thus a conservative approximation of \eqref{eq:prob_cc_distrob_true_ind}.
    The condition on $N$ in the theorem's statement ensures that $\delta<\alpha$. \hfill $\qed$ 

\end{document}